\documentclass[twoside,a4paper,10pt]{article}
\usepackage{latexsym,amssymb}
\usepackage{enumerate}
\usepackage{graphicx}
\usepackage{color}
\usepackage{hyperref}
\usepackage[frenchb,english]{babel}
\usepackage[matrix,arrow]{xy}
\def \be{\begin{eqnarray*}}
\def \ee{\end{eqnarray*}}
\def \ben{\begin{enumerate}}
\def \een{\end{enumerate}}
\def \beit{\begin{itemize}}
\def \eeit{\end{itemize}}
\def \bui#1#2{\mathrel{\mathop{\kern 0pt#1}\limits^{#2}}}
\def \buil#1#2{\mathrel{\mathop{\kern 0pt#1}\limits_{#2}}}
\def \bfll{\begin{flushleft}}
\def \efll{\end{flushleft}}
\def \bflr{\begin{flushright}}
\def \eflr{\end{flushright}}

\def \findemo{\hfill$\square$\\}

\def \Span{\mathrm{span}}

\def \SO{\mathrm{SO}}

\def \U{\mathrm{U}}

\def \R{\mathbb{R}}

\def \C{\mathbb{C}}
\def \sgn{\mathrm{sgn}}
\def \Ric{\mathrm{Ric}}
\def \S{\mathrm{S}}

\def \ph{\varphi}

\def \lrc{\lrcorner}
\def \clN{\cdot_{\scriptscriptstyle N}}
\newcommand\ip{{\langle\cdot \,,\cdot \rangle}}

\newcommand{\la}{\langle}       
\newcommand{\ra}{\rangle}

\newtheorem{ethm}{Theorem}[section]
\newtheorem{edefi}[ethm]{Def\mbox{}inition}
\newtheorem{elemme}[ethm]{Lemma}
\newtheorem{erem}[ethm]{Note}
\newtheorem{re}[ethm]{Remark}

\newtheorem{ecor}[ethm]{Corollary}
\newtheorem{prop}[ethm]{Proposition}
\newtheorem{eexemple}[ethm]{Example}

\title{Skew Killing spinors in four dimensions}
\author{Nicolas Ginoux\footnote{Universit\'e de Lorraine, CNRS, IECL, F-57000 
Metz, France,  E-mail: \texttt{nicolas.ginoux@univ-lorraine.fr}},\, Georges 
Habib\footnote{Lebanese University, Faculty of Sciences II, Department of 
Mathematics, P.O. Box 90656 Fanar-Matn, Lebanon,
E-mail: \texttt{ghabib@ul.edu.lb}},\, Ines Kath\footnote{Universit\"at 
Greifswald, Institut f\"ur Mathematik und Informatik, 
Walther-Rathenau-Stra\ss{}e 47 17487 Greifswald, Germany, E-mail: 
\texttt{ines.kath@uni-greifswald.de}
}}

\begin{document}
\maketitle
\noindent\begin{center}\begin{tabular}{p{115mm}}
\begin{small}{\bf Abstract.} This paper is devoted to the classification of 
$4$-dimensional Riemannian spin manifolds carrying skew Killing spinors. A skew 
Killing spinor $\psi$ is a spinor that satisfies the equation 
$\nabla_X\psi=AX\cdot\psi$ with a skew-symmetric endomorphism $A$. 
We consider the degenerate case, where the rank of $A$ is at most two everywhere 
and the non-degenerate case, where the rank of $A$ is four everywhere. We prove 
that in the degenerate case the manifold is locally isometric to the Riemannian 
product $\R\times N$ with $N$ having a skew Killing spinor and we explain under 
which conditions on the spinor  the special case of a local isometry to 
$\mathbb{S}^2\times\mathbb{R}^2$ occurs. In the non-degenerate case, the 
existence of skew Killing spinors is related to doubly warped products whose 
defining data we will describe.  
\end{small}\\
\end{tabular}\end{center}

\noindent\begin{small}{\it Mathematics Subject Classification} (2010): 53C25, 
53C27.  
\end{small}

\noindent\begin{small}{\it Keywords}: Generalized Killing spinors, doubly warped 
product, Hodge operator. 
\end{small}

\section{Introduction}\label{s:intro}

Let $(M^n,g)$ be an $n$-dimensional Riemannian spin manifold.
A generalised Killing spinor on $M$ is a section $\psi$ of the spinor bundle 
$\Sigma M$ of $M$ satisfying the overdetermined differential equation 
$\nabla_X\psi=AX\cdot\psi$
for some symmetric endomorphism field $A$ of $TM$.
Here and as usual, ``$\cdot$'' denotes the Clifford multiplication on $\Sigma 
M$.
Numerous papers have been devoted to the classification of Riemannian spin 
manifolds carrying such spinors.
Several results have been obtained for particular $A$ but it is still an open 
problem to get a complete classification for general $A$.
Let us quote some of these results. 
First, recall that when $A$ is the zero tensor field, that is, the corresponding 
spinor is parallel, then McK. Wang  \cite{Wang89}  showed that such manifolds 
can be characterised by their holonomy groups which can be read off the Berger 
classification.
The case where $A$ is a nonzero real multiple of the identity is that of 
classical real Killing spinors. It was shown by C.~B\"ar \cite{Baer93} that real 
Killing spinors correspond to parallel spinors on the (irreducible) cone over 
the manifold, to which then McK.\,Wang's result applies. Furthermore, in 
dimension $n\le 8$, there are several results on a classification up to isometry 
\cite{BFGK,Hijazi86}.
When the tensor $A$ is parallel \cite{Morel03}, or a Codazzi tensor \cite{BGM05} 
or both $A$ and $g$ are analytic \cite{BMM13} (see also \cite{CS07}), it is 
shown that the manifold  $M$
is isometrically embedded into another spin manifold of dimension $n+1$ carrying 
a parallel spinor and that the tensor $A$ is the half of the second fundamental 
form of the immersion.
We also cite the partial classification of generalised Killing spinors on the 
round sphere \cite{MorSemmsphere14, MorSemmspherethree14} and on $4$-dimensional 
Einstein manifolds of positive scalar curvature \cite{MorSemmEinstein12} where 
in some cases the generalised Killing spinor turns out to be a Killing spinor.

In this paper, we are interested in an equation dual to the generalised Killing 
one, which we call {\it skew Killing spinor equation}.
More precisely, on a given Riemannian spin manifold $(M^n,g)$, a spinor field 
$\psi$ is called a skew Killing spinor if it satisfies for some {\it 
skew-symmetric} endomorphism field $A$ of $TM$ the differential equation
\begin{equation}\label{eq:sks}
\nabla_X\psi=AX\cdot\psi
\end{equation}
for all $X\in TM$. This equation was originally defined in \cite{HabibRoth12}. 
 Each skew Killing spinor is a parallel section with respect to the modified 
metric connection $\nabla-A\otimes\mathrm{Id}$, in particular it has constant 
length. 
Moreover, for a given skew symmetric endomorphism field $A$ of $TM$, 
the space of skew Killing spinors is a complex vector space of dimension at 
most $\mathrm{rk}_{\mathbb{C}}(\Sigma M)=2^{\left[n/2\right]}$.

Very few examples of Riemannian spin manifolds $(M^n,g)$ carrying skew Killing 
spinors are known for which $A\neq0$.
For $2$-dimensional manifolds, apart from $\R^2$ or quotients thereof with 
trivial spin structure, only the round sphere of constant curvature can carry 
such spinors and in that case they correspond to restrictions of Killing spinors 
from $\mathbb{S}^3$ onto totally geodesic $\mathbb{S}^2$  \cite{HabibRoth12}.
In that case, the tensor $A$ coincides with the standard complex structure $J$ 
induced by the conformal class of $\mathbb{S}^2$ or with $-J$ depending on the 
sign of the Killing constant chosen on $\mathbb{S}^3$. Each skew Killing spinor 
on $\mathbb{S}^2$ immediately gives rise to a three-dimensional example, namely 
to a skew Killing spinor on $\mathbb{S}^2\times \R$, where $A=\pm J$ on 
$\mathbb{S}^2$ is trivially extended to the $\R$-factor.
More generally, for a manifold of dimension $n=3$ the following is known 
\cite[Prop. 4.3]{HabibRoth12}.  
If $M^3$ admits a skew Killing spinor~$\psi$, 
then, locally, $\psi$ can be transformed into a parallel spinor by a suitable 
conformal change of the metric. 
In particular, $M^3$ is locally conformally 
flat. 
If, in addition, $M^3$ is simply-connected, then this conformal change is 
defined globally.  
Conversely, if $(M^3,g)$ admits a nonzero parallel spinor, 
then for any conformal change of $g$, there exists a skew Killing spinor with 
respect to the new metric. 
See Section~\ref{Exex} for more detailed information.

In dimensions $6$ and $7$, there are lots of examples provided by 
$\mathrm{SU}(3)$- resp. $\mathrm{G}_2$-structures on $M$, see e.g. types 
$\chi_1,\chi_2,\chi_4$ in  \cite[Lemma 
3.5]{ACFH15} and type $\mathcal{W}_2,\mathcal{W}_4$ in \cite[Lemma 4.5]{ACFH15} 
respectively.
\color{black}

Obvious examples in four dimensions can be obtained as products $N\times \R$, 
where $N$ is a three-dimensional manifold admitting a skew Killing spinor, see 
Example~\ref{ex1}. A special case of this construction is the product  
$\mathbb{S}^2\times\mathbb{R}^2$, see Example~\ref{ex2}. For each of the 
endomorphisms $A^\pm:=\pm J\oplus{0}$, this manifold admits the maximal number 
of skew Killing spinors. 

The main purpose of this work is to establish a classification result when the 
dimension of $M$ is four. Note that the pointwise rank of $A$ is either zero, 
two or four.
We will split the classification into two parts. In Section~\ref{S4} we will 
study the degenerate case, where the rank of $A$ is at most two everywhere.
In Section~\ref{S5} we will consider the case where $\mathrm{rk}(A)=4$ on all of 
$M$. Before we start the classification, we determine the general integrability 
conditions in arbitrary dimensions arising from the existence of a skew Killing 
spinor, see Section \ref{s:integrcond}.
In Sections \ref{S3} and \ref{S4}, we specify these 
conditions to four dimensions, especially to the degenerate case. We use that 
the spinor bundle $\Sigma M$ splits into the eigenspaces $\Sigma^+M$ and 
$\Sigma^-M$ of the volume form and the bundle of two-forms splits into those of 
self-dual and of anti-self-dual forms, which act on $\Sigma^\pm M$. 
We also adapt some techniques used in \cite{MorSemmEinstein12} but for a 
skew-symmetric endomorphism $A$.
We use the integrability conditions to achieve the following classification 
result in case that the Killing map is degenerate everywhere. 

{\bf Theorem A.} {\it
\ Let $(M^4,g)$ be a connected Riemannian spin manifold carrying a skew Killing 
spinor~$\psi$, where the rank of the corresponding skew-symmetric tensor field 
$A$ is at most two everywhere.
Then either $\psi$ is parallel on $M$ or, around every point of $M$, we have a 
local Riemannian splitting $\R\times N$ with $N$ having a skew Killing spinor. 
If, in addition, the length of the summand $\psi^+$ in the decomposition 
$\psi=\psi^++\psi^-\in\Sigma^+M\oplus\Sigma^-M$ is not constant, then we are in 
the second case with $N=\R\times \mathbb{S}^2$, that is, $(M,g)$ is a local 
Riemannian product $\mathbb{S}^2\times \R^2$ around every point. } 

For a more detailed formulation see Theorem~\ref{finaldeg}, where we also 
discuss the global structure of $(M,g)$ if $M$ is complete. 

Let us turn to the case where the Killing map is non-degenerate everywhere. In 
Section~\ref{S51} we will prove that, essentially, the existence of a skew 
Killing spinor $\psi$ with non-degenerate Killing map $A$ is equivalent to the 
existence of a Killing vector field $\eta$ and an almost complex structure $J$ 
satisfying certain conditions, see Proposition~\ref{P1} for a detailed 
formulation. The spinor $\psi$ and the data $\eta$ and $J$ are related by the 
equations $J(X)\cdot\psi^-=iX\cdot\psi^-$ and $g(\eta,X)=\langle X\cdot 
\psi^+,\psi^-\rangle/|\psi|^2$ for all $X\in TM$.

In Section~\ref{S52}, we consider the special case where $A\eta$ is parallel to 
$J\eta$. Then $AJ=JA$ holds and $J$ is integrable, see Remark~\ref{rint}. 
Manifolds with skew Killing spinors satisfying these conditions are related to 
doubly warped products. A \emph{doubly warped product} is a Riemannian manifold 
$(M,g)$ of the form $(I\times \hat 
M,dt^2\oplus\rho(t)^2\hat{g}_{\hat{\eta}}\oplus\sigma(t)^2\hat{g}_{\hat{\eta}
^\perp})$, where  $(\hat M,\hat{g})$ is a Riemannian manifold with unit Killing 
vector field $\hat{\eta}$, and $\hat{g}_{\hat{\eta}}$, 
$\hat{g}_{\hat{\eta}^\perp}$ are the components of the metric $\hat{g}$ along 
$\R\hat{\eta}$ and $\hat{\eta}^\perp$, respectively, $I\subset\R$ is an open 
interval and $\rho,\sigma\colon I\to\R$ are smooth positive functions on $I$. 
Locally, doubly warped products can be equivalently described as local 
DWP-structures, see the appendix. On $\hat M$, we define a function $\hat\tau$ 
by $\hat \nabla_X\hat\eta=\hat\tau\cdot \hat J(X)$ for $X\in\hat \eta^\perp$, 
where $\hat J$ is a fixed Hermitian structure on $\eta^\perp$. Locally, $(\hat 
M,\hat g)$ is a Riemannian submersion over a two-dimensional base manifold $B$. 
Let $\hat K$ denote the Gaussian curvature of $B$.  We obtain the following 
result, see Theorem~\ref{thm} and Corollary~\ref{corfinal}.
 
{\bf Theorem B.} {\it
Let $(M,g)$ admit a skew Killing spinor such that  $A\eta || J\eta$ and 
$|\eta|\not\in\{0,1/2\}$ everywhere. Then $M$ is locally isometric to a doubly 
warped product for which the data $\hat K$ and $\hat \tau$ are constant and 
$\rho$ and $\sigma$ satisfy the differential equations
$$
(\sigma^2)'=-\frac 2{\sqrt{1-4\rho^2}}\rho\hat \tau,\qquad 
(\sigma^2)'\,\frac{\rho'}\rho=\hat K-2\frac{\rho^2}{\sigma^2}\hat\tau^2.
$$
Conversely, if $M$ is isometric to a simply-connected doubly warped product for 
which the data $\hat K$ and $\hat \tau$ are constant and $\rho$ and $\sigma$ 
satisfy the above differential equations, then $(M,g)$ admits a skew Killing 
spinor such that $A\eta || J\eta$.
}

The differential equations in Theorem B can be locally solved and one obtains 
explicit formulas for the doubly warped product. Let us finally mention that the 
skew Killing spinors on $M=I\times \hat M$ are related to quasi Killing spinors 
in the sense of \cite{FriedKim00} on $\hat M$, see Remark~\ref{notequasi}.

{\bf Acknowledgement:} The second named author would like to thank the 
Alexander 
von Humboldt foundation and the DAAD for the financial support.

\section{General integrability conditions for skew Killing 
spinors}\label{s:integrcond}

In this section we give a few necessary conditions for the existence of nonzero 
skew Killing spinors. Before we state the main result, we recall some facts from 
Riemannian and spin geometry, see e.g. \cite[Chap. 1]{BFGK} or \cite[Chap. 
2]{LM}.

In all this paper we identify, on a Riemannian manifold $(M^n,g)$, one-forms 
with vector fields via the metric $g$.
Recall that the Hodge star operator is defined by
\[\omega\wedge *\omega'=\langle\omega,\omega'\rangle\mathrm{vol}_g\]
for all differential $p$-forms $\omega,\omega'$ on $M$, where $\mathrm{vol}_g$ 
is the volume form of $M$ (giving its orientation). The Hodge star operator 
satisfies $*^2=(-1)^{p(n-p)}$ on $p$-forms and has the following useful 
properties 
\begin{equation}\label{l:basics}
X\wedge *\omega=(-1)^{p+1}*(X\lrcorner\,\omega) \quad\textrm{and}\quad 
X\lrcorner *\omega=(-1)^p *(X\wedge\omega)
\end{equation}
for any vector field $X$. Recall also that the Clifford multiplication between a 
vector field $X$ and a differential $p$-form $\omega$ is defined as  
\begin{equation}\label{basics2}
X\cdot\omega=X\wedge \omega-X\lrcorner\,\omega \quad\textrm{and}\quad 
\omega\cdot X=\omega\wedge X+(-1)^p X\lrcorner\, \omega,
\end{equation}
from which the identity $X\cdot Y\cdot +Y\cdot X\cdot=-2g(X,Y)$ follows for any 
vector fields $X$ and $Y.$ 

From now on, we assume $M$ to be spin with fixed spin structure.
In that case, there exists a Hermitian vector bundle $\Sigma M\to M$, called the 
spinor bundle, on which the tangent bundle $TM$ acts by Clifford multiplication, 
$TM\otimes \Sigma M\to \Sigma M; X\otimes\psi\mapsto X\cdot\psi.$ We will write 
$XY\cdot\psi$ instead of $X\cdot Y\cdot\psi$. 
Recall that a real $p$-form also acts by Clifford multiplication in a formally 
self- or skew-adjoint way according to its degree: for any $p$-form $\omega$ and 
any spinors $\varphi,\psi$, we have 
\[\langle\omega\cdot\varphi,\psi\rangle=(-1)^{\frac{p(p+1)}{2}}
\cdot\langle\varphi,\omega\cdot\psi\rangle.\] 

The Levi-Civita connection $\nabla$ on $M$ defines a metric connection, also 
denoted by $\nabla$, on $\Sigma M$ with respect to the Hermitian product $\ip$ 
and that preserves Clifford multiplication.
In other words, for all $X,Y\in\Gamma(TM),$ the rules
$$X(\langle 
\psi,\varphi\rangle)=\langle\nabla_X\psi,\varphi\rangle+\langle\psi,
\nabla_X\varphi\rangle,\qquad 
\nabla_X(Y\cdot\varphi)=\nabla_XY\cdot\varphi+Y\cdot\nabla_X\varphi$$
are satisfied for all spinor fields $\psi,\varphi$.
If we denote by $R_{X,Y}:=[\nabla_X,\nabla_Y]-\nabla_{[X,Y]}$ the curvature 
tensor associated with the connection $\nabla$, the spinorial Ricci identity 
states that, for all $\psi$ and $X$,
\begin{equation}\label{eq:spinricci}
-\frac{1}{2}\mathrm{Ric}(X)\cdot\psi=\sum_{j=1}^n e_j\cdot 
R_{X,e_j}\psi,\end{equation}
see e.g. \cite[Eq. 1.13]{BFGK}. 

In the following, we will assume the manifold $M$ to carry a skew Killing 
spinor 
field $\psi$ with corresponding skew-symmetric endomorphism $A$.
We make $A$ into a $2$-form via the metric $g$, that is, we consider 
$(X,Y)\mapsto g(AX,Y)$, which we still denote by $A$. 
In a pointwise orthonormal basis $\{e_i\}_{i=1,\cdots,n}$ of $TM$, we have 
$A=\frac{1}{2}\sum_{j=1}^n e_j\wedge Ae_j$ (mind the factor $\frac{1}{2}$).  
In particular, Clifford multiplication of any spinor field $\psi$ by $A$ is 
given by
\begin{equation}\label{e:2form}
A\cdot \psi={\textstyle \frac12}\sum_{j=1}^n e_j \cdot Ae_j\cdot\psi \,.
\end{equation}

In the next proposition, we compute the curvature data arising from the 
existence of such a spinor.
These integrability equations will play a crucial role for the classification 
in 
the $4$-dimensional case.

\begin{prop}\label{p:integrcond}
Let $\psi$ be any solution of {\rm (\ref{eq:sks})} on a spin manifold $(M^n,g)$ 
for some skew-symmetric endomorphism field $A$ of $TM$.
Then the following identities hold for $X,Y\in \Gamma(TM)$
\begin{enumerate}
\item $R_{X,Y}\psi=\left((\nabla_XA)(Y)-(\nabla_Y A)(X)+2AY\wedge 
AX\right)\cdot\psi$.
\item $-\frac{1}{2}\mathrm{Ric}(X)\cdot\psi=\left(\nabla_XA+X\lrcorner\, 
dA+(\delta A)(X)+4A\wedge AX+2A^2X\right)\cdot\psi$, where $d$ is the exterior 
derivative and $\delta$ is the codifferential w.r.t. the metric $g$.
\item $\mathrm{S}\cdot\psi=4\left(2dA+\delta A+4A\wedge 
A+|A|^2\right)\cdot\psi$, where $\S$ denotes the scalar curvature of $(M,g)$ 
and 
$|A|^2:=\sum_{j=1}^n|Ae_j|^2$ written in any pointwise orthonormal basis 
$(e_j)_{1\leq j\leq n}$ of $TM$.
\end{enumerate}
\end{prop}

{\it Proof}: We derive (\ref{eq:sks}) and take suitable traces of the 
identities 
obtained. 
First, if $x\in M$ and $X,Y\in\Gamma(TM)$ such that $\nabla X=\nabla Y=0$ at 
$x$, then
\begin{eqnarray*}
\nabla_X\nabla_Y\psi&=&\nabla_X(AY\cdot\psi)
\ =\ (\nabla_XA)(Y)\cdot\psi+AY\cdot\nabla_X\psi\\
&=&(\nabla_XA)(Y)\cdot\psi+AY\cdot AX\cdot\psi
\end{eqnarray*}
at $x$. Thus, with the help of Equations (\ref{basics2}), we write
\begin{eqnarray*} 
R_{X,Y}\psi&=&\nabla_X\nabla_Y\psi-\nabla_Y\nabla_X\psi\\
&=&\big((\nabla_XA)(Y)-(\nabla_Y A)(X)+AY\cdot AX-AX\cdot AY\big)\cdot\psi\\
&=&\big((\nabla_XA)(Y)-(\nabla_Y A)(X)+2AY\wedge 
AX-g(AY,AX)+g(AX,AY)\big)\cdot\psi\\
&=&\big((\nabla_XA)(Y)-(\nabla_Y A)(X)+2AY\wedge AX\big)\cdot\psi,
\end{eqnarray*}
which is the first identity.

Next we fix a local orthonormal basis of $TM$, which we denote by $(e_j)_{1\leq 
j\leq n}$.
Using the spinorial Ricci formula (\ref{eq:spinricci}) and the identities 
(\ref{basics2}), we compute
\begin{eqnarray*}
-\frac{1}{2}\mathrm{Ric}(X)\cdot\psi&=&\sum_{j=1}^n e_j\cdot R_{X,e_j}\psi 
\ =\ \sum_{j=1}^ne_j\cdot\left((\nabla_XA)(e_j)-(\nabla_{e_j} A)(X)+2Ae_j\wedge 
AX\right)\cdot\psi\\
&=&\Big(\sum_{j=1}^ne_j\cdot(\nabla_XA)(e_j)-\sum_{j=1}^ne_j\wedge 
(\nabla_{e_j} 
A)(X)+\sum_{j=1}^ne_j\lrcorner(\nabla_{e_j} A)(X)\\
&&\phantom{\Big\{}+2\sum_{j=1}^ne_j\cdot (Ae_j\wedge AX)\Big)\cdot\psi.
\end{eqnarray*}
Now we compute each term separately. 
First, 
$
\sum_{j=1}^ne_j\cdot (\nabla_XA)(e_j)\cdot\psi=2\nabla_XA\cdot \psi
$
by (\ref{e:2form}), where we see $\nabla_XA$ as a $2$-form on $M$.
The second sum can be computed in terms of the exterior and the covariant 
derivatives of $A$. Namely
\begin{eqnarray*}
-\sum_{j=1}^ne_j\wedge (\nabla_{e_j} A)(X)&=&\Big(\sum_{j=1}^n e_j\wedge 
\nabla_{e_j}A\Big)(X)-\sum_{j=1}^ng(X,e_j)\nabla_{e_j}A
\ =\ X\lrcorner\, dA-\nabla_XA.
\end{eqnarray*}
The third sum can be expressed in terms of the codifferential of $A$:
\[\sum_{j=1}^ne_j\lrcorner(\nabla_{e_j} A)(X)=\sum_{j=1}^n(\nabla_{e_j} 
A)(X,e_j)=-\sum_{j=1}^n(\nabla_{e_j} A)(e_j,X)=(\delta A)(X).\]
It remains to notice that, by Equations (\ref{basics2}), we have 
\begin{eqnarray*}
\sum_{j=1}^ne_j\cdot (Ae_j\wedge AX)\cdot\psi&=&\sum_{j=1}^n e_j\cdot Ae_j\cdot 
AX\cdot\psi + \sum_{j=1}^n {g(Ae_j,AX)} e_j\cdot\psi \\
&=&(2 A\cdot AX -A^2X)\cdot\psi\ =\ (2 A\wedge AX+A^2X )\cdot\psi.
\end{eqnarray*}
This shows the second equation.

To obtain the scalar curvature, we trace the spinorial Ricci identity. Given a 
local orthonormal basis $(e_j)_{1\leq j\leq n}$ of $TM$, we write
\begin{eqnarray*}
\frac{\mathrm{S}}{2}\psi&=&-\frac{1}{2}\sum_{j=1}^n e_j\cdot 
\mathrm{Ric}(e_j)\cdot\psi\\
&=&\sum_{j=1}^n e_j\cdot\big(\nabla_{e_j}A+e_j\lrcorner\, dA+(\delta 
A)(e_j)+4A\wedge Ae_j+2A^2e_j\big)\cdot\psi\\
&\bui{=}{(\ref{basics2})}&\sum_{j=1}^n 
\Big(e_j\wedge\nabla_{e_j}A-e_j\lrcorner\nabla_{e_j}A\Big)\cdot\psi+\sum_{j=1}
^n\Big(e_j\wedge(e_j\lrcorner\, dA)-\underbrace{e_j\lrcorner(e_j\lrcorner\, 
dA)}_{0}\Big)\cdot\psi\\
&&+\sum_{j=1}^n(\delta A)(e_j) e_j\cdot\psi+4\sum_{j=1}^n \Big( e_j\wedge 
A\wedge Ae_j-\underbrace{e_j\lrcorner(A\wedge Ae_j)}_{0}\Big)\cdot\psi\\&&
+2\sum_{j=1}^n \Big( \underbrace{e_j\wedge 
A^2e_j}_{0}-g(A^2e_j,e_j)\Big)\cdot\psi\\
&=&\left(dA+\delta A+3dA+\delta A+8A\wedge A+2|A|^2\right)\cdot\psi\\
&=&\left(4dA+2\delta A+8A\wedge A+2|A|^2\right)\cdot\psi,
\end{eqnarray*}
which is the last identity.
Here, we use the the identity $\sum_{j=1}^n e_j\wedge (e_j\lrcorner\, 
\omega)=p\omega$, which holds for any $p$-form $\omega$.
\findemo

\section{The vector fields $\eta$ and $\xi$ in four dimensions}\label{S3}

In this section, we consider a $4$-dimensional spin manifold $(M,g)$ that 
carries a skew Killing spinor.
On spin manifolds of even dimension $2m$, the complex volume form 
$(\mathrm{vol}_g)_\C:=i^me_1\cdot e_2\ldots\cdot e_{2m}$, where 
$(e_j)_{j=1,\cdots,2m}$ is an arbitrary orthonormal frame, splits the spinor 
bundle into two orthogonal subbundles that correspond to the eigenvalues $\pm 
1$ 
of $(\mathrm{vol}_g)_\C.$
Hence, on our four-dimensional manifold $(M,g)$, we have $\Sigma M=\Sigma^+M 
\oplus\Sigma^-M $, where 
$$\Sigma^\pm M:= \{\psi\in \Sigma M\mid (\mathrm{vol}_g)_\C\cdot\psi=\pm 
\psi\}. 
$$
The spaces $\Sigma^\pm M$ are preserved by the connection $\nabla$ of the 
spinor 
bundle and are interchanged by Clifford multiplication by tangent vectors.
According to this decomposition, we write any spinor field $\psi$ as 
$\psi=\psi^++\psi^-$ and we set $\bar\psi:=\psi^+-\psi^-$. 
Recall now that differential forms act on the spinor bundle $\Sigma M$ as 
follows: for any differential $p$-form $\omega$ on $M$ and 
$\psi\in\Gamma(\Sigma 
M)$ 
\begin{equation}\label{l:b1}
\omega\cdot\psi=*\omega\cdot \bar\psi\quad {\rm for}\quad p=1,2\quad {\rm 
and}\quad  \omega\cdot\psi=- (* \omega) \cdot\bar\psi\quad {\rm for}\quad p=3,4.
\end{equation} 
Let $\bigwedge_{\pm}^2M=\{\omega\in\bigwedge^2M\mid *\omega=\pm\omega\}$ be the 
spaces of self-dual and anti-self-dual forms on $M$.
For $\omega\in\bigwedge^2M$ we denote by $\omega_\pm$ the projections of 
$\omega$ to these spaces.
Then, one can easily see from Equations (\ref{l:b1}) that $\bigwedge_{\pm}^2M$ 
acts trivially on $\Sigma^{\mp}M$ and that the maps
\begin{equation}\label{action}\textstyle
\bigwedge^2_-M \longrightarrow \Sigma^-M\cap (\psi^-)^\perp,\ 
\omega_-\longmapsto \omega_-\cdot \psi^-, \quad \bigwedge^2_+M  \longrightarrow 
\Sigma^+M\cap (\psi^+)^\perp,\ \omega_+\longmapsto \omega_+\cdot \psi^+
\end{equation}
are isomorphisms if $\psi^+\not=0$ and $\psi^-\not=0$.

Now assume that $\psi$ is a skew Killing spinor of norm one. By decomposing 
$\psi$ into $\psi^+$ and $\psi^-$ as we said before, 
we obtain isomorphisms (\ref{action}) on the open set $M':=M_0\cap M_1,$ with 
$$M_0:=\{x\in M\mid \psi^-(x)\not=0\}\quad {\rm and} \quad M_1:=\{x\in M\mid 
\psi^+(x)\not=0\}.$$

Equation (\ref{eq:sks}) can be written as $\nabla_X\psi^\pm=AX\cdot\psi^\mp.$  
We define a  vector field $\eta$ on $M$ and a vector field $\xi$ on $M_0$ by 
\begin{equation}\label{xieta}
g(\eta,X):=\langle X\cdot \psi^+,\psi^-\rangle,\quad \psi^+=:\xi\cdot\psi^-,
\end{equation}
where the definition of $\xi$ uses that the map $T_pM\to \Sigma^+_pM$, 
$X\mapsto 
X\cdot\psi^-$ is bijective at each~$p\in M_0.$
Then, clearly $\eta=-|\psi^-|^2\xi$ holds on $M_0$ and  
$1=|\psi^+|^2+|\psi^-|^2=|\psi^-|^2(1+|\xi|^2)$. We define
$$f:=1-2|\psi^-|^2, \quad \rho:=|\eta|\le 1/2.$$
Then
\begin{equation}
\rho=\frac{|\xi|}{1+|\xi|^2},\quad f=\frac{|\xi|^2-1}{|\xi|^2+1},\quad 
f^2=1-4\rho^2,\quad \eta=\frac12(f-1)\xi \label{rhoxi}
\end{equation}
holds, where these functions are defined. 

We collect some properties of $\eta$ and $\xi$ that will be used later on. 
\begin{elemme}\label{l:h} On $M$, we have
\begin{enumerate}
\item $df=4 A\eta $ \label{i:eta1}
\item $\nabla_X\eta = f AX$\label{i:eta2}
\item $d\eta=2fA,\,\,\,\delta\eta=0$ \label{i:eta3}
\item $fdA=-4A\eta \wedge A\,.$\label{i:eta4}
\end{enumerate}
\end{elemme}

{\it Proof}: Differentiating the function $|\psi^-|^2$ along any vector field 
$X\in TM$ gives
$$X(|\psi^-|^2)=2\la\nabla_X\psi^-,\psi^-\ra=2\la AX\cdot 
\psi^+,\psi^-\ra=2g(\eta,AX)=-2g(A\eta,X).$$
This proves {\it \ref{i:eta1}}.
To prove {\it \ref{i:eta2}}, we consider two vector fields $X$ and $Y$ that can 
be assumed to be parallel at some point $x\in M$ to compute
\begin{eqnarray*} 
g(\nabla_X \eta,Y) &=& X(g(\eta,Y))\ =\ X(\langle Y\cdot \psi^+,\psi^-\rangle)\\
&=& \la Y\cdot AX\cdot \psi^-,\psi^-\ra +\la Y\cdot \psi^+,AX\cdot \psi^+\ra\\
&=& -g(Y,AX)|\psi^-|^2+g(Y,AX)|\psi^+|^2\\
&=&(1-2|\psi^-|^2) g(AX,Y)
\end{eqnarray*}
at $x$, which is {\it \ref{i:eta2}}.
Moreover,
$$d\eta(X,Y)=(\nabla_X\eta)(Y)-(\nabla_Y\eta)(X)=2fA(X,Y),$$
which yields the first part of {\it \ref{i:eta3}}.
The divergence of $\eta$ is clearly zero by {\it \ref{i:eta2}} and the fact 
that 
$A$ is skew-symmetric. Finally, 
$$ 0=dd\eta=2fdA+2df\wedge A,$$
which together with {\it\ref{i:eta1}} gives {\it \ref{i:eta4}}. \findemo
\begin{re} \label{rxi}It follows from Lemma \ref{l:h} that $\nabla\eta$ is 
skew-symmetric on $M$ which means that $\eta$ is a Killing vector field on $M$. 
\end{re}

The open sets $M_0$ and $M_1$ are dense in $\{p\in M\mid A_p\not=0\}$. Indeed, 
if, e.g., $\psi^-$ vanishes on some open set $U\subset \{p\in M\mid 
A_p\not=0\}$, then so does its covariant derivative and therefore 
$AX\cdot\psi^+=0$ on $U$. Hence $A=0$ on $U$, which contradicts the assumption 
on $A$.  

With the notation introduced above, we have 
$M_0=\{x\in M\mid f(x)\not=1\}$ and $M_1=\{x\in M\mid f(x)\not=-1\}$. Then 
$M'=M_0\cap M_1=\{x\in M\mid f(x)\not=\pm1\}=\{x\in M\mid \rho(x)\not=0\}$. We 
define also the set
$$M'':=\left\{x\in M\mid \rho(x)\not\in\{0,\textstyle{\frac12}\}\right\}=M'\cap 
\left\{x\in M\mid \rho(x)\not=\textstyle{\frac12}\right\} =\{x\in M\mid 
f(x)\not\in\{0,\pm1\} \}.$$
By Lemma~\ref{l:h}, 1., the open set $M''$ is dense in $\{p\in M\mid 
A_p(\eta)\not=0\}$. In particular, $M''\subset M$ is dense if $A$ is 
non-degenerate everywhere. The case where $\rho=1/2$ on an open set will be 
treated in  Proposition~\ref{h=1/2}.

\begin{re}\label{33} Let us change the orientation of $M$ and denote by 
$\hat\Sigma M$ the spinor bundle with respect to the new orientation. Then we 
can identify $\hat\Sigma M$ with $\Sigma M$ via $\hat\Sigma^+M=\Sigma^-M$ and 
$\hat\Sigma^+M=\Sigma^-M$ Accordingly, we define a section $\hat\psi$ of 
$\hat\Sigma M$ by $\hat\psi^+=\psi^-$, $\hat\psi^-=\psi^+$. With $\psi$ also 
$\hat \psi$ is a skew Killing spinor and  the vector fields $\hat \xi$ and 
$\hat 
\eta$ associated with $\hat\psi$ are equal to $\hat\xi=-\xi/|\xi|^2$ and 
$\hat\eta=-\eta$, respectively. 

On $M''$, we have $|\xi|\not=1$. Hence, if there exists a skew Killing spinor 
on 
$M$ and if $M=M''$, then we always may assume that $|\xi|>1$ up to a possible 
change of orientation on each connected component of $M$. If  $|\xi|>1$, then 
$f$ is positive, thus $f=\sqrt{1-4\rho^2}$. 
\end{re}

\section{The degenerate case}\label{S4}
In this section, we assume that $\mathrm{rk}(A)\leq2$ everywhere on $M^4$, 
which 
is equivalent to suppose that the kernel of $A$ is at every point either $4$- 
or 
$2$-dimensional.  Then $AX\wedge A=0$ for all $X\in TM$. In particular, $dA=0$ 
on $M''$ by Lemma~\ref{l:h}. 

\subsection{Examples}\label{Exex}

\begin{eexemple}\label{ex1}
If $N$ is a $3$-dimensional spin manifold with a skew Killing spinor $\varphi$, 
then $N\times \mathbb{R}$ admits a skew Killing spinor $\psi\not=0$ for which 
$|\psi^+|=|\psi^-|$~holds.
\end{eexemple}
Let us prove the above statement. Recall that the spinor bundle of $M=N\times 
\mathbb{R}$ is given by $\Sigma M=\Sigma N\oplus \Sigma N$ and the Clifford 
multiplication on $M$ is related to the one on $N$ by \cite{Baer98} 
$$(X\clN\oplus -X\clN)\psi=X\cdot\partial_t\cdot\psi.$$   
where $\partial_t$ is the unit vector field on $\mathbb{R}$ and $X\in TN$. Now 
we set $\psi:=\varphi+\partial_t\cdot\varphi$ according to the above 
decomposition. Let $A$ denote the Killing map associated with $\psi$. Then we 
can easily check that $\nabla_{\partial t}\psi=0$ and, for $X\in TN$, 
\begin{eqnarray*}
\nabla_X\psi&=&\nabla_X\varphi+\partial_t\cdot\nabla_X\varphi\\
&=&AX\clN\varphi+\partial_t\cdot (AX\clN\varphi)\\
&=&AX\cdot\partial_t\cdot\varphi+\partial_t\cdot AX\cdot\partial_t\cdot\varphi\\
&=& AX\cdot\psi.
\end{eqnarray*}
Hence $\psi$ is a skew Killing spinor on $M$. The vector field $\xi$ in this 
example is just $-\partial_t$ which is parallel. Since $|\partial_t|=1$, we 
have 
$|\psi^+|=|\psi^-|$.

Let us recall at this point, what is known about three-dimensional manifolds 
with skew Killing spinors. As already mentioned in the introduction, each skew 
Killing spinor on $\mathbb{S}^2$ immediately gives rise to a three-dimensional 
example, namely to a skew Killing spinor on $\mathbb{S}^2\times \R$. 
Furthermore, if $\dim N=3$  and  if $(N,g)$ admits a skew Killing 
spinor~$\psi$, 
then $N$ is locally conformally flat \cite[Prop. 4.3]{HabibRoth12}. Indeed, 
locally, there exists a function $u$ such that $\psi$ transforms into a 
parallel 
spinor $\bar \psi$ with respect to the metric $\bar g :=e^{2u}g$ and 
three-dimensional Riemannian manifolds with a non-trivial parallel spinor field 
are flat. If $N$ is simply-connected, then $u$ is globally defined. In the 
latter case the metric $\bar g$ is not necessarily complete even if $(N,g)$ is. 
 
Conversely, if $(N,g)$ admits a nonzero parallel spinor, then for any conformal 
change of the metric on the manifold $N$ there exists a skew Killing spinor 
with 
respect to the new metric. We conclude this overview with the flat case 
$N=\R^3$. If $\psi\not=0$ is a solution of {\rm (\ref{eq:sks})} on 
$N=\mathbb{R}^3$ endowed with the flat metric, then $A=0$ and $\psi$ is a 
parallel spinor field. Indeed, as mentioned above, there exists a globally 
defined function $u$ on $\R^3$ such that the metric $\bar g :=e^{2u}g$ admits a 
parallel spinor. Hence, $\bar g$ is also flat. In particular, the scalar 
curvature $\bar \S$ vanishes. On the other hand, $\bar \S= 
8e^{-2u}e^{-u/2}\Delta e^{u/2}$ since $\bar g$ arises by conformal change from 
the flat metric $g$. Thus  $\Delta(e^{u/2})=0$, that is, $e^{u/2}$ is a 
harmonic 
function on $\mathbb{R}^3$.
But since $e^{u/2}\geq0$, Liouville's theorem implies that $e^{u/2}$ -- and so 
$u$ itself -- is constant.
This shows~$A=0$.

\begin{eexemple}\label{ex2}
We consider $M=\mathbb{S}^2\times\mathbb{R}^2$. Let $J$ denote the standard 
complex structure on $\mathbb{S}^2$. We define endomorphisms $A^{\pm}:=\pm 
J\oplus{0}$ on $TM=T\mathbb{S}^2\oplus T\mathbb{R}^2$. For each of these 
endomorphisms, the space of skew Killing spinors is four-dimensional. It can be 
spanned by elements with non-vanishing $A\eta$ and it also can be spanned by 
elements for which $A\eta=0$ holds. 
\end{eexemple}

Let us prove this statement. The spinor bundle of 
$\mathbb{S}^2\times\mathbb{R}^2$ is pointwise given by 
$\Sigma(\mathbb{S}^2\times\mathbb{R}^2)=\Sigma\mathbb{S}^2\otimes 
\Sigma\mathbb{R}^2$ and the Clifford multiplication on 
$\mathbb{S}^2\times\mathbb{R}^2$ is \cite{Baer98} 
$$X\cdot(\varphi\otimes\sigma)=(X\cdot_{\mathbb{S}^2}\varphi)\otimes\bar\sigma,
\qquad Y\cdot(\varphi\otimes\sigma)=\varphi\otimes 
(Y\cdot_{\mathbb{R}^2}\sigma), $$
for $X\in T\mathbb{S}^2$ and $Y\in T\mathbb{R}^2.$
Now, we consider on $\mathbb{S}^2$ a skew Killing spinor $\varphi$, 
corresponding to the standard complex structure $J$, and a parallel spinor 
$\sigma$ in $\Sigma^+(\mathbb{R}^2)$ of norm $1$.
The spinor field $\psi:=\varphi\otimes\sigma$ is clearly a skew Killing spinor, 
since in the $\mathbb{S}^2$-direction we have 
$$\nabla_X\psi=(\nabla_X\varphi)\otimes\sigma=(JX\cdot_{\mathbb{S}^2}
\varphi)\otimes\sigma=JX\cdot(\varphi\otimes\sigma)=JX\cdot\psi$$
and $\nabla_Y\psi=0$ in the $\mathbb{R}^2$-direction.
The same computation holds when replacing $J$ by $-J$ and choosing $\sigma\in 
\Sigma^-(\mathbb{R}^2)$.
As the spaces of skew Killing spinors $\varphi$ corresponding to the standard 
complex structure $J$ or its opposite on $\mathbb{S}^2$ are each complex 
$2$-dimensional, we deduce that the space of skew Killing spinors with Killing 
map $A^+$ is at least -- and therefore exactly -- $4$-dimensional. The same 
holds for $A^-$.
In particular, each skew Killing spinor on $\mathbb{S}^2\times\mathbb{R}^2$ is a 
linear combination with constant coefficients of skew Killing spinors for $A^+$ 
and also one of skew Killing spinors for $A^-$. 
Note that the vector field $\xi$, associated to the above-defined skew Killing 
spinor $\psi$, is the one coming from the spinor $\varphi$ on $\mathbb{S}^2,$ 
since $T\mathbb{S}^2\simeq\Sigma^+\mathbb{S}^2$ and
$$\psi^+=\varphi^+\otimes\sigma=(\xi_{\mathbb{S}^2}
\cdot\varphi^-)\otimes\sigma=\xi_{\mathbb{S}^2}
\cdot(\varphi^-\otimes\sigma)=\xi_{\mathbb{S}^2}\cdot\psi^-.$$
Therefore, $\xi=\xi_{\mathbb S^2}$ and $A^2\xi=J^2\xi_{\mathbb 
S^2}=-\xi_{\mathbb S^2}$, which cannot vanish on the sphere. Thus 
$A\eta\not=0$. 
If we consider instead of the above constructed $\psi$ the spinor $\psi+Y\cdot 
\bar \psi$
for a parallel vector field $Y$ on $\R^2$ with $|Y|=1$, we obtain a skew Killing 
spinor with $\xi=-Y$, hence $A\eta=0$. 

\subsection{Classification}
Let us first assume that $\rho=1/2$ on an open set. By definition of $\rho$, 
this condition is equivalent to $|\psi^+|=|\psi^-|$. 
We prove that, under this assumption, the manifold is locally isometric to that 
in Example \ref{ex1}.
\begin{prop}\label{h=1/2}
Let $\psi$ be a nonzero skew Killing spinor on $M^4$ and assume that 
$|\psi^+|=|\psi^-|$ on an open set $U$. Then $U$ is a local Riemannian product 
of a line by a $3$-dimensional Riemannian manifold carrying a skew Killing 
spinor.    
\end{prop}

{\it Proof}: Let $\psi$  be a skew Killing spinor of norm one such that 
$|\psi^+|=|\psi^-|$. Then $f=0$ by definition of $f$.  Thus $\eta$ is parallel 
by Lemma~\ref{l:h}. In this case $\eta^\perp$ is integrable and the spinor 
$\psi$ restricts to a skew Killing spinor on the integral manifolds.  In fact, 
for any given integral manifold $N$, its spinor bundle is identified with 
$\Sigma^+ M$, so the spinor $\varphi=\psi^+$ restricts to a skew Killing spinor 
on~$N$. Indeed,
$$\nabla^N_X\varphi=\nabla_X^M\psi^+=AX\cdot\psi^-=-AX \cdot \xi \cdot 
\psi^+=-AX\clN\varphi,
$$ which proves the assertion. \findemo

In the next part of the section, we want to exclude the case $\rho=1/2$ and make 
the stronger assumption

\hspace{2cm}\parbox{10cm}{$(M^4,g)$ is a Riemannian spin manifold carrying a 
skew Killing spinor such that $M=M''$ and $\mathrm{rk}(A)=2$ everywhere.}    
\hfill{(GA)}

Due to the orthogonal splitting of the spinor bundle $\Sigma 
M=\Sigma^+M\oplus\Sigma^-M$ we can decompose further the equations in 
Proposition \ref{p:integrcond} in order to get more integrability conditions. 
Namely, 

\begin{elemme} \label{pro:intecon} Under the assumption (GA), we have
\begin{eqnarray}
0&=&\textstyle{\frac12} \Ric(X)+2A^2X+*(\xi\wedge \nabla_X A)+\xi\lrcorner 
\nabla _X A +(\delta A)(X)\xi \label{e:1}\\[1ex]
0&=&\Big( \textstyle{\frac12} \Ric(X)\wedge \xi  +2A^2X\wedge \xi +\nabla_X 
A\Big)_- \label{e:2}\\[1ex]
0&=&\textstyle{\frac12} \Ric(\xi)+2A^2\xi-\delta A\label{e:3}\\
0&=&\delta A+(|A|^2-\textstyle{\frac14}{\rm S} )\xi\label{e:4}\\
0&=&\Big( \xi\wedge\delta A\Big)_-\label{e:5}\\
0&=&-(\delta A)(\xi)+|A|^2-\textstyle{\frac14} {\rm S} \label{e:6}
\end{eqnarray}
 for any vector field $X$.
\end{elemme}
{\it Proof}: We take the orthogonal projection of the formulas in Proposition 
\ref{p:integrcond}  to $\Sigma^+M$ and $\Sigma^-M$.
This gives, after using $\psi^+=\xi\cdot\psi^-,$  $dA=0$ and $A\wedge AX=0$ that
\begin{eqnarray}
0&=&\left(\textstyle{\frac12} \Ric(X) +2 A^2X\right)\cdot\psi^- + \left 
(\nabla_X A +(\delta A)(X)\right)\xi\cdot\psi^- \label{e:I+}\\
0&=&\left(\textstyle{\frac12} \Ric(X) +2 A^2X\right)\cdot\xi\cdot\psi^- + \left 
(\nabla_X A +(\delta A)(X)\right)\cdot\psi^-\label{e:I-}
\end{eqnarray}
and
\begin{eqnarray}
0&=&\left(|A|^2-\textstyle{\frac14} {\rm S}\right)\cdot\xi\cdot\psi^- + (\delta 
A)\cdot \psi^-\label{e:II+}\\
0&=&\left(|A|^2-\textstyle{\frac14} {\rm S} \right)\cdot\psi^- + (\delta 
A)\cdot\xi\cdot \psi^-,\label{e:II-}
\end{eqnarray}
respectively.
Equation (\ref{e:I+}) gives 
\begin{eqnarray*}
0&=&\left(\textstyle{\frac12} \Ric(X) +2 A^2X+\xi\wedge \nabla_X A+\xi\lrc 
\nabla_XA  +(\delta A)(X)\xi   \right)\cdot \psi^-.
\end{eqnarray*}
Hence, by formula (\ref{l:b1}), we obtain (\ref{e:1}).
Equation (\ref{e:I-}) yields 
\begin{eqnarray*}
0&=&\left(\textstyle{\frac12} \Ric(X)\wedge \xi -\textstyle{\frac12} \Ric(X,\xi) 
+2 A^2X\wedge\xi-2g(A^2X,\xi)+\nabla_X A +(\delta A)(X)  \right)\cdot \psi^-.
\end{eqnarray*} 
Now, by taking the scalar product with $\psi^-$ and identifying the real part, 
the $0$-th order term must vanish.
This is Equation (\ref{e:3}).
Also, we have
$$ \left( \textstyle{\frac12} \Ric(X)\wedge \xi  +2A^2X\wedge \xi +\nabla_X 
A\right)\cdot\psi^-=0.$$
The isomorphism from $\bigwedge ^2_-M$ to the orthogonal complement 
$(\psi^-)^\perp$ yields Equality (\ref{e:2}) from the above identity.
Equation (\ref{e:II+})  gives (\ref{e:4}). Finally, Equation (\ref{e:II-}) 
yields
$$0=\left(|A|^2-\textstyle{\frac14} {\rm S} -(\delta A)(\xi)-\xi\wedge \delta 
A\right)\cdot \psi^-.$$
Taking the Hermitian product with $\psi^-$, we obtain Equations (\ref{e:5}) and 
(\ref{e:6}) after identifying the real parts.
\findemo

In the following, we will further simplify the equations in Lemma 
\ref{pro:intecon}.

\begin{prop}\label{prop:intcond2}
Under the assumption (GA), we have
\begin{eqnarray}
\delta A&=& 0\label{eq:deltaA}\\
\S&=& 4|A|^2 \label{e:6'}\\
 \Ric(\eta)&=&-4A^2\eta\label{eq:Ric}\\
\nabla_\eta A&=&0\label{eq:nablaxi}\\
 (\nabla_XA)(\eta)&=&-f\left(\textstyle{\frac14}\Ric(X) +A^2X\right) 
\label{eq:xinabla}\\
\label{eq:nablaXAxi}
\nabla_X(A\eta)&=&-\frac{f}{4}{\rm Ric}(X),\\
\label{eq:nablastara}
\eta\lrcorner \nabla_X (*A)&=&\ \nabla_X((*A)\eta)\ =\ 
\frac{1}{4}\mathrm{Ric}(X)+A^2X
\end{eqnarray}
for every $X\in TM$.
\end{prop}

{\it Proof}: By (\ref{e:5}), we have
$$\xi\wedge \delta A=*(\xi \wedge \delta A)=-\xi\lrc *\delta A.$$
Hence, the interior product with $\xi$  yields $0=\xi\lrc (\xi\wedge \delta 
A)=|\xi|^2\delta A-(\delta A)(\xi)\cdot \xi.$
 Now, applying Equation (\ref{e:4}) to $\xi$ gives 
$$0=(\delta A)(\xi)+\left(|A|^2-\textstyle{\frac14} {\rm S}\right)|\xi|^2,$$
which, after combining with (\ref{e:6}), leads to
$0=(1+|\xi|^2)(\delta A)(\xi),$
which gives (\ref{eq:deltaA}). Now (\ref{e:4}) yields (\ref{e:6'}).
Equation (\ref{eq:Ric}) now follows from (\ref{e:3}) and (\ref{eq:deltaA}). 

From (\ref{e:2}), we get
$$
*\Big( \textstyle{\frac12} \Ric(X)\wedge \xi  +2A^2X\wedge \xi +\nabla_X A\Big)
=\textstyle{\frac12} \Ric(X)\wedge \xi  +2A^2X\wedge \xi +\nabla_X A,
$$
which, by Equation (\ref{l:basics}), is equivalent to 
$$
\textstyle{\frac12} \xi\lrc *(\Ric(X)) +2\xi\lrc*(A^2X) +*\nabla_X A
=\textstyle{\frac12} \Ric(X)\wedge \xi  +2A^2X\wedge \xi +\nabla_X A.
$$
Taking the interior product by $\xi$, this gives
$$\xi\lrc*\nabla_X A
= \xi\lrc\left(\left(\textstyle{\frac12}\Ric(X)+2A^2X\right)\wedge \xi\right) 
+\xi\lrc \nabla_X A,$$
thus
$$*(\xi\wedge\nabla_X A)=\left(-\textstyle{\frac12}\Ric(X)-2A^2X\right)|\xi|^2  
+\xi\lrc \nabla_X A$$
by Equations (\ref{l:basics}) and (\ref{eq:Ric}). On the other hand, Equations 
(\ref{e:1}) and (\ref{eq:deltaA}) give
$$*(\xi\wedge \nabla_X A)=-\textstyle{\frac12} \Ric(X) -2A^2X-\xi\lrcorner 
\nabla _X A.$$
Substracting and adding the latter two equations and replacing $\xi$ by 
$-(1+|\xi|^2)\eta$ yields (\ref{eq:xinabla}) and the identity $\eta\lrcorner 
\nabla_X (*A)=\frac{1}{4}\mathrm{Ric}(X)+A^2X$ for all $X\in \Gamma(TM)$. The 
last equation yields (\ref{eq:nablastara}) since 
$(*A)(\nabla_X\eta)=f(*A)(AX)=*fAX\wedge A=0$.
Furthermore, Equation (\ref{eq:xinabla}) shows that the expression $(\nabla_X 
A)(\eta,Y)$ is symmetric in $X$ and $Y$.
Thus
$$0=(\nabla_X A)(\eta,Y)-(\nabla_Y A)(\eta,X)
=-dA(X,Y,\eta)+(\nabla_\eta A)(X,Y)=(\nabla_\eta A)(X,Y)$$
by $dA=0$. 
This proves (\ref{eq:nablaxi}). 
Equation (\ref{eq:nablaXAxi}) follows from (\ref{eq:xinabla}) together with 
$\nabla _X\eta=fAX$.
\findemo

\begin{re}
We can prove integrability conditions analogous to those in 
Lemma~\ref{pro:intecon} and Proposition~\ref{prop:intcond2} also for arbitrary 
rank of $A$. These general conditions are more involved. Since we will not use 
them in the present paper, we do not state them here.
\end{re}

\begin{elemme}\label{p:dim4Axi=0}
Under the assumption (GA), the set $\{p\in M\mid A\eta|_p\not=0\}$ is dense in 
$M$.
\end{elemme}

{\it Proof:} 
Assume that $A\eta=0$ on an open set $U$. We know that $\eta$ is a Killing 
vector field on $M$.
Moreover, by Lemma \ref{l:h}, the vector field $\eta$ has constant length on 
$U$. Indeed, for every $X\in TM$,
\[X\left(|\eta|^2\right)=2g\left(\nabla_X\eta,\eta\right)=2fg(AX,
\eta)=-2fg(A\eta,X)=0.\]
By \cite[Thm. 4]{BeresNiko08}, since (\ref{eq:Ric}) implies $\Ric(\eta)=0$, we 
can conclude that $\eta$ is parallel on $U$. But this contradicts item 2 of 
Lemma~\ref{l:h} since $f\not=0$ and $A\not=0$ everywhere by assumption.
\findemo

In the following, we will often assume assume that $A\eta\neq0$ on all of $M$. 
If $A\eta\neq0$, then we have $A^2\eta\neq0$ everywhere, thus the vectors 
$\frac{A\eta}{\left|A\eta\right|}$ and $\frac{A^2\eta}{\left|A^2\eta\right|}$ 
form an orthonormal basis of the image of $A$.
As $A$ is of rank $2$, we obtain
\begin{equation}\label{eq:AintermsofAxiA2xi} 
A=\frac{1}{|A\eta|^2}\cdot A\eta\wedge A^2\eta.
\end{equation}
Furthermore,  note that (\ref{eq:AintermsofAxiA2xi}) already implies
 \begin{equation}\label{eq:a3}
 A^3\eta=-\frac{|A^2\eta|^2}{|A\eta|^2}A\eta=-\frac \S 8 A\eta,
 \end{equation}
where the last equality comes from the identity (\ref{e:6'}). Obviously, 
$A^3\eta=-\frac \S 8 A\eta$ holds also if $A\eta=0$. 

Since  $df=4A\eta$ by Lemma \ref{l:h}, (\ref{eq:nablaXAxi}) implies
\begin{equation}\label{eq:nabladf=-fric}
\nabla df=-f\cdot\Ric.
\end{equation}
This equation has been extensively studied in \cite{GinouxHabibKath18}. 
Using this formula, we now express the Ricci tensor of the vector field 
$A\eta$.

\begin{elemme}\label{p:Nflat}  If (GA) holds, then the Ricci tensor satisfies  
\begin{equation}\label{3.5}
\mathrm{Ric}(A\eta)=\frac{\mathrm{S}}{2}\cdot A\eta+\frac{f}{16}\cdot 
d\mathrm{S},\qquad {\rm Ric}((*A)\eta)=\frac{1}{16}d\S.
\end{equation}
In particular, we have
\begin{equation}\label{Svergl}
(A\eta)({\rm S})=f((*A)\eta)({\rm S}).
\end{equation}
\end{elemme}

{\it Proof:} By Bochner's formula for $1$-forms, $\Delta(df)-{\rm 
Ric}(df)=\nabla^*\nabla (df)$ holds. Since $\nabla df=-f{\rm Ric}$ is symmetric 
and since $\nabla^*=\delta$ on symmetric $(0,2)$-tensors, this gives
$$
\Delta(df)-{\rm Ric}(df)=\delta\nabla df=\delta(-f{\rm Ric})=\Ric (df)-f\delta 
({\rm Ric})
=\Ric(df) +\frac f2d{\rm S},
$$
where we used the well-known identity $d\S=-2\delta \Ric .$
Hence, we deduce
$$
\Delta(df)=2\mathrm{Ric}(df)+\frac{f}{2}d{\rm S}.
$$
But (\ref{eq:nabladf=-fric}) also gives $\Delta f=-\mathrm{tr}_g\left(\nabla 
df\right)=f\cdot\mathrm{S}$, so that $\Delta(df)=d(\Delta 
f)=d(f\cdot\mathrm{S})$. Therefore
$$
\mathrm{Ric}(df)=\frac{1}{2}\mathrm{S}\cdot df+\frac{f}{4}\cdot d\mathrm{S}.
$$
The first equation in (\ref{3.5}) now follows from the equality $df=4A\eta$.

In the following, we will compute the Ricci curvature of the vector field 
$(*A)\eta.$ Notice first that $(*A)\eta=\eta\lrcorner(*A)=*(\eta\wedge A)$. 
Hence, this vector field belongs to the kernel of $A$ as  
$$g(AX,(*A)\eta){\rm vol}_g=AX\wedge *^2(\eta\wedge A)=-AX\wedge \eta\wedge 
A=0$$
for any $X\in TM$. Based on the fact $A\eta\lrcorner(*A)=*(A\eta\wedge A)=0$, we 
first compute
\begin{equation}\label{rest}
A\eta\lrc \nabla_ X(*A)= -(*A)(\nabla_X A\eta) =\frac f4 (*A)(\Ric(X))=*\frac f4 
(\Ric(X) \wedge A).
\end{equation}
This gives
$$
\eta\lrc(A\eta \lrc\nabla_X(*A)) = -*\frac{f}{4}(\eta\wedge {\rm Ric}(X)\wedge 
A) = -\frac{f}{4}{\rm Ric}(X)\lrcorner*(\eta \wedge A)
 = -\frac{f}{4}{\rm Ric}((*A)\eta,X).
$$
On the other hand, by (\ref{eq:nablastara}) and (\ref{eq:a3}), we have
$$
A\eta\lrc(\eta\lrcorner\nabla_X (*A))  =  A\eta\lrcorner \big(\frac{1}{4}{\rm 
Ric}( X)+A^2X\big)
 = g\big(\frac{1}{4}{\rm Ric}(A\eta)+A^3\eta,X\big)=\frac{f}{64}g(d{\rm S},X).
$$
Comparing the two identities gives the second equation in (\ref{3.5}).
Equation (\ref{Svergl}) can be deduced from computing ${\rm 
Ric}(A\eta,(*A)\eta))$ in two ways from (\ref{3.5}) taking the scalar product by 
$(*A)\eta$ in the first formula and by $A\eta$ in the second one. Remember that 
$(*A)\eta$ lies in the kernel of $A$.
\findemo

In the following, we will establish and prove three technical lemmas (Lemmas 
\ref{l:nablaaeta}, \ref{lem:riccitensoratwoeta} and \ref{lem:equationscal}), 
which will show that the kernel and the image of the endomorphism $A$ are 
integrable and totally geodesic. Then the proof of Theorem~A  will follow from 
the de Rham theorem. 

\begin{elemme}\label{l:nablaaeta}  Assume that (GA) holds.
Then we have the identity 
\begin{equation}
\nabla_{A\eta}A^2\eta=-\frac{f}{4}{\rm Ric}(A^2\eta)-\frac{f^2}{32}A(d{\rm S}). 
\label{pp}
\end{equation}
\end{elemme}
{\it Proof:} 
By continuity, it suffices to prove the assertion on the set $\{p\in M\mid 
A\eta|_p\not=0\}$ since this set is dense in $M$ by Lemma \ref{p:dim4Axi=0}. 
Thus we may assume that $A\eta\not=0$ everywhere. For any $X\in TM$, we have 
\begin{equation}\label{daeta}
d(|A\eta|^2)=2g(\nabla (A\eta),A\eta)=-\frac f2 \Ric(A\eta),
\end{equation}
where we use Equation (\ref{eq:nablaXAxi}) in the last equality.
Thus, from Lemma \ref{p:Nflat}, we find  
$$d\left(\frac{1}{|A\eta|^2}\right)=\frac f{2|A\eta|^4}\Ric(A\eta)=\frac 
f{4|A\eta|^4}\big( \S\, \cdot A\eta+\frac f8 d\S\big). $$
Moreover, $\delta(A^2\eta)=0$. Indeed, for any two-form $\omega$ in four 
dimensions and any vector $X$, the formula 
$\delta(X\lrc\omega)=*(dX\wedge*\omega)-\delta\omega(X)$ holds. Using $\delta 
A=0$ and $4d(A\eta)=ddf=0$, this yields
\begin{eqnarray*}\label{deltaatwo}
\delta(A^2\eta)=\delta(A\eta\lrc A)=*(d(A\eta)\wedge*A)-(\delta A)(A\eta)=0.
\end{eqnarray*}
Now, by taking the divergence of both sides of (\ref{eq:AintermsofAxiA2xi}), we 
compute 
\begin{eqnarray*}
0 &=& \delta A\ =\ \delta\left(\frac1{|A\eta|^2} A\eta\wedge A^2\eta\right)\ =\ 
-d\left(\frac1{|A\eta|^2}\right) \lrc (A\eta\wedge A^2\eta) + \frac1{|A\eta|^2} 
\delta(A\eta\wedge A^2\eta)\\
&=& -\frac f{4|A\eta|^2}(\S\, A\eta+\frac f8 d\S)\lrc A +\frac 1{|A\eta|^2}\big( 
\delta(A\eta) A^2\eta+\nabla_{A^2\eta} A\eta -\nabla_{A\eta}{A^2 
\eta}-\delta(A^2\eta) A\eta\big),
\end{eqnarray*}
where we use the formula $\delta(X\wedge Y)=(\delta X)Y+\nabla_Y X-\nabla_X 
Y-(\delta Y)X,$ valid for any $X,Y\in TM$. Furthermore, the divergence of 
$A\eta$ is equal to $f\S/4$ as an easy consequence from tracing Equation 
(\ref{eq:nablaXAxi}). This finally gives (\ref{pp}).
\findemo

The following technical lemma expresses a partial trace of the Ricci tensor.
\begin{elemme} \label{lem:riccitensoratwoeta} Assume that (GA) holds and that 
$A\eta\not=0$ everywhere. 
Then the following identity holds:
$$\frac1{|A\eta|^2}{\rm Ric}(A\eta,A\eta)+\frac1{|A^2\eta|^2}{\rm 
Ric}(A^2\eta,A^2\eta)=\S-\frac2{f\S}A\eta({\rm S}).$$
\end{elemme}
{\it Proof:} The proof relies on taking the scalar product of ${\rm 
Ric}(A^2\eta)$ in Lemma \ref{l:nablaaeta} with the vector field $A^2\eta.$
Indeed, we have 
\begin{eqnarray*} 
{\rm Ric}(A^2\eta,A^2\eta)&=&-\frac 
4f\Big(g(\nabla_{A\eta}A^2\eta+\frac{f^2}{32}A(d\S),A^2\eta)\Big)\\
&=& -\frac2f A\eta(|A^2\eta|^2)+\frac f8 g(d\S, A^3\eta)\\
&\bui{=}{(\ref{eq:a3})}& -\frac 2f A\eta \Big(\frac \S 8 
|A\eta|^2\Big)-\frac{f\S}{64}A\eta(\S)\\
& \bui{=}{(\ref{daeta})} &\frac 
\S8\Ric(A\eta,A\eta)-\Big(\frac{|A\eta|^2}{4f}+\frac{f\S}{64}\Big)A\eta(\S). 
\end{eqnarray*}
Hence, again by (\ref{eq:a3}), we find 
$$\frac{\Ric(A^2\eta,A^2\eta)}{|A^2\eta|^2}=\frac{\Ric(A\eta,A\eta)}{|A\eta|^2}
-\left(\frac2{f\S}+\frac f{8|A\eta|^2}\right)A\eta(\S).$$
Finally, the identity
$$\frac{\Ric(A\eta,A\eta)}{|A\eta|^2}=\frac\S2+\frac f{16|A\eta|^2}A\eta(\S),$$
which follows from Lemma \ref{p:Nflat}, leads to the required equality.
\findemo

\begin{elemme} \label{lem:equationscal}  If (GA) holds, then the scalar 
curvature is constant and ${\rm Ric}+4A^2=0.$
\end{elemme}
{\it Proof:}  As in the proof of Lemma \ref{l:nablaaeta}, we may assume that 
$A\eta\not=0$ everywhere.  By Lemma \ref{p:Nflat}  we know that 
$${\rm Ric}(A\eta) -\frac{\rm S}2 A\eta = f\cdot {\rm Ric}((*A)\eta).$$ We take 
the divergence of both sides. We start with the left hand side. Note that for 
any vector field $X\in \Gamma(TM)$ the formula $\delta({\rm 
Ric}(X))=g(\delta{\rm Ric},X)-\sum_{i=1}^n g({\rm Ric}(e_i),\nabla_{e_i}X)$ 
holds, where $e_1,\dots,e_n$ is any pointwise orthonormal basis. Using this and 
$\delta(A\eta)=\frac{f{\rm S}}{4}$, we compute
\begin{eqnarray}
\delta({\rm Ric}(A\eta)-\frac{{\rm S}}{2}A\eta)&=&g(\delta{\rm 
Ric},A\eta)-\sum_{i=1}^4 g({\rm Ric}(e_i),\nabla_{e_i}A\eta)- 
\frac{1}{2}\big(-g(d{\rm S}, A\eta)+{\rm S}\,\delta(A\eta)\big) \nonumber\\
&=&-\frac{1}{2}g(d{\rm S},A\eta)+\frac{f}{4}\sum_{i=1}^4 g({\rm Ric}(e_i),{\rm 
Ric}(e_i))+\frac{1}{2}g(d{\rm S}, A\eta)-\frac{f}{8}{\rm S}^2\nonumber\\
&=&\frac{f}{4}|{\rm Ric}|^2-\frac{f}{8}{\rm S}^2.\label{eq:divergencericciaeta}
\end{eqnarray}
To get the divergence of the right hand side, we first compute that of the 
vector field ${\rm Ric}((*A)\eta)$.
For this, we use the same formula as above and again $d\S=-2\delta \Ric$ to 
write
\begin{eqnarray}
\delta({\rm Ric}((*A)\eta))&=&-\frac{1}{2}((*A)\eta)({\rm S})-\sum_{i=1}^4 
g({\rm Ric}(e_i),\nabla_{e_i}((*A)\eta))\nonumber\\
&=&-\frac{1}{2f}(A\eta)({\rm S})-\sum_{i=1}^4 g({\rm 
Ric}(e_i),\nabla_{e_i}((*A)\eta)). \label{xx}
\end{eqnarray}
In the last equality, we used (\ref{Svergl}). Inserting (\ref{eq:nablastara}) 
into (\ref{xx}), we find  
$$
\delta({\rm Ric}((*A)\eta))=-\frac{1}{2f}(A\eta)({\rm S})-\frac{1}{4}|{\rm 
Ric}|^2-\sum_{i=1}^4 g({\rm Ric}(e_i),A^2e_i), 
$$
which in turn gives
\begin{eqnarray}
\delta(f\cdot {\rm Ric}((*A)\eta))&=&-g(df,{\rm Ric}((*A)\eta)) +f\cdot \delta 
({\rm Ric}((*A)\eta)))\nonumber\\
& =& \ -\frac{3}{4}(A\eta)({\rm S})-\frac{f}{4}|{\rm Ric}|^2-f\sum_{i=1}^4 
g({\rm Ric}(e_i),A^2e_i) \label{eq:deltascal2}
\end{eqnarray}
by (\ref{3.5}). Comparing Equations (\ref{eq:divergencericciaeta}) and 
(\ref{eq:deltascal2}), we obtain 
$$\sum_{i=1}^4 g({\rm Ric}(e_i),A^2e_i)=-\frac{3}{4f}(A\eta)({\rm 
S})-\frac{1}{2}|{\rm Ric}|^2+\frac{1}{8}{\rm S}^2. $$
On the other hand, this sum can be computed on the particular orthonormal frame 
$\frac{A\eta}{|A\eta|},\frac{A^2\eta}{|A^2\eta|},e_3, e_4$ with $e_3,e_4$ in the 
kernel of $A$ as follows: using Lemma \ref{lem:riccitensoratwoeta}, we write
\begin{eqnarray}
\sum_{i=1}^4 g({\rm Ric}(e_i),A^2e_i)&=&\frac{1}{|A\eta|^2}{\rm 
Ric}(A\eta,A^3\eta)+\frac{1}{|A^2\eta|^2}{\rm Ric}(A^2\eta,A^4\eta)\nonumber\\
&\bui{=}{(\ref{eq:a3})}&-\frac \S8\left(\frac{1}{|A\eta|^2}{\rm 
Ric}(A\eta,A\eta)+\frac{1}{|A^2\eta|^2}{\rm 
Ric}(A^2\eta,A^2\eta)\right)\nonumber\\
&=&-\frac {\S^2}8 + \frac1{4f}A\eta({\rm S}).\label{fast1}
\end{eqnarray}
Comparing these two computations yields
\begin{equation}\label{fast2}
4 A\eta({\rm S}) =f ({\rm S}^2-2 |{\rm Ric}|^2).
\end{equation}
The Cauchy-Schwarz Inequality gives 
\begin{equation}
\sum_{i=1}^4 g({\rm Ric}(e_i),A^2e_i)\le |\Ric||A^2|. \label{CSU}
\end{equation} 
We take the square of this inequality. Then we use (\ref{fast1}) and 
(\ref{fast2}) to express the left and the right hand side, respectively. We 
obtain 
$$\Big(-\frac {\S^2}8 + \frac1{4f}A\eta({\rm 
S})\Big)^2\le\Big(\frac{\S^2}2-\frac2f A\eta(\S)\Big)\cdot\frac{\S^2}{32}= 
\frac{\S^4}{64} -\frac{\S^2}{16f}A\eta(\S),$$
where besides (\ref{e:6'}), which says that $\S= 4|A|^2$,  we used 
$|A^2|^2=(|A|^2)^2/2$, which follows from the fact that $A$ is skew-symmetric of 
rank two. This inequality is only true if $A\eta(\S)=0$. But then (\ref{CSU}) is 
an equality. Hence, $\Ric$ is a multiple of $A^2$ at every point of $M^4$.  
Since Tr\,$\Ric=\S$ and Tr$A^2=-|A|^2=-\S/4$, we obtain $\Ric=-4A^2.$ As the 
vector field $(*A)\eta$ lies in the kernel of $A$, the second equation in 
(\ref{3.5}) implies that the scalar curvature is constant. This ends the proof.
\findemo

\begin{elemme}\label{loc}
If (GA) is satisfied, then $(M,g)$ is locally isometric to $\R^2\times 
\mathbb{S}^2$.
\end{elemme}
{\it Proof:}  We show that the two orthogonal distributions ${\rm Im}(A)$ and 
${\rm Ker}(A)$ -- which are both of rank two by assumption -- are parallel. If 
this is proved to be true, then we get a local Riemannian product by the de Rham 
decomposition theorem. Clearly, it suffices to show that ${\rm Im}(A)$ is 
parallel since ${\rm Ker}(A)={\rm Im}(A)^\perp$.  Let us first consider the open 
subset $V:=\{p\in M\mid A\eta|_p\not=0\}$. On $V$, the image of $A$ is spanned 
by $A\eta$ and $A^2\eta$. Note that  $\nabla_X A\eta=fA^2X$ by 
(\ref{eq:nablaXAxi}) and Lemma~\ref{lem:equationscal}. Thus $\nabla_XA\eta$ is 
contained in ${\rm Im}(A)$ for all $X\in TM$.
Furthermore, by Equation~(\ref{rest}) and Lemma~\ref{lem:equationscal}, we have  
$A\eta\lrc \nabla_ X(*A)=0$ for all $X\in TM$. Equation~(\ref{e:2}) now gives 
$A\eta\lrc \nabla_ XA=0$. Thus  
$\nabla_XA^2\eta= A(\nabla_X(A\eta))=fA^3(X)$. In particular, also 
$\nabla_XA^2\eta$ is contained in ${\rm Im}(A)$ for all $X\in TM^4$. This proves 
that ${\rm Im}(A)$ is parallel.

We want to extend this splitting of $TM$ into two parallel distributions to all 
of $M$. To this end, we observe that, on $V$, the Ricci map has constant 
eigenvalues $0,0,\S/2,\S/2>0$ and Ker$(A)$ and Im$(A)$ are the 
eigendistributions. Since $V\subset M$ is dense by Lemma~\ref{p:dim4Axi=0}, 
these are also the eigenvalues of $\Ric$ on all of $M$ and the two-dimensional 
eigendistributions of $\Ric$ are parallel on all of $M$.
We deduce that $(M,g)$ is locally isometric to the Riemannian product 
$\mathbb{R}^2\times\mathbb{S}^2$.
\hfill$\square$

Now we can prove the main result of this section. In particular, it says that, 
in the degenerate case, the skew Killing spinor is parallel or $(M,g)$ is 
locally isometric to one of the examples discussed in Section~\ref{Exex}.

\begin{ethm} \label{finaldeg}
Let $(M^4,g)$ be a connected Riemannian spin manifold carrying a skew Killing 
spinor $\psi$, where the rank of the corresponding skew-symmetric tensor field 
$A$ is $\le 2$ everywhere.
Then either $\psi$ is parallel (i.e., $A=0$) on $M$ or, around every point of 
$M$, we have a local Riemannian splitting $\R\times N$ with $N$ having a skew 
Killing spinor.  
If, moreover,  $|\psi^+|$ (thus also $|\psi^-|$) is not constant, then $(M,g)$ 
is a local Riemannian product $\mathbb{S}^2\times \R^2$ around every point and 
the Killing map  equals~$\pm J\oplus {\rm 0}$.  

If, in addition, $(M,g)$ is complete, then $(M,g)$ is globally isometric to the 
Riemannian product $\mathbb{S}^2\times\Sigma^2$, where $\Sigma^2$ is either flat 
$\R^2$, a flat cylinder with trivial spin structure or a flat $2$-torus with 
trivial spin structure. 
\end{ethm}

{\it Proof:} We define $U:=\{p\mid A_p\not=0\}$ and $U':=U\cap M'$, $U'':=U\cap 
M''$. Recall that $U'\subset U$ is dense. We know that 
Equation~(\ref{eq:nabladf=-fric}) holds on the open set $U''$. We claim that it 
holds on all of $M$. Obviously, it is true on the closure $\overline{U''}$ of 
$U''$. It also holds on $U'\setminus  \overline{U''}$ since this set is open 
with $f\equiv 0$. Consequently, it holds on $U'$, thus on $U$ since $U'\subset 
U$ dense.  Hence it is true on supp$(A)=\overline U$. Furthermore, on the 
complement of supp$(A)$, we have $df=0$ and $\Ric=0$,  thus 
(\ref{eq:nabladf=-fric}) holds on $M$.  Now we can apply Prop.\,1.2 in 
\cite{HePetersenWylie15}, which shows that either $f\equiv 0$ on $M$ or 
supp$(f)=M$. If $f\equiv 0$, then Proposition~\ref{h=1/2} applies. Assume now 
that supp$(f)=M$. Then $M''$ is dense in $M$.  Let $U$ and $U''$ be defined as 
above.  On $U''$, the assumption (GA) is satisfied. As we have seen, the 
eigenvalues of $\Ric$ are $0$ and $\S/2$ and the eigendistributions of $\Ric$ on 
$U''$ are parallel. Thus this holds also on $\overline U''=\overline U$.
If $\overline U=M$, then we are done by Lemma~\ref{loc}. If $\overline 
U=\emptyset$, then $\psi$ is parallel. Assume that $\overline U$ were non-empty 
and not equal to $M$. Then the complement $W$ of $\overline U$ is open and not 
empty with $A=0$. Thus $\psi$ is parallel on $W$, hence $\Ric=0$ on $W$, thus 
also on $\overline W$. Since $M$ is connected, $\overline U\cap \overline W$  is 
non-empty. Hence we can chose a point $p$ in this intersection. But then 
$p\in\overline U$ would imply that $\S/2>0$ is an eigenvalue of $\Ric_p$ and 
$p\in\overline W$ would imply that $\Ric_p=0$, a contradiction.

Note that, as we already noticed in \cite[Theorem 2.4]{GinouxHabibKath18}, the 
manifold $(M,g)$ must be globally isometric to the product 
$\mathbb{S}^2\times\Sigma^2$, where $\Sigma^2$ is a quotient of flat $\R^2$.
The reason is that the fundamental group of $M$ can act on the 
$\mathbb{S}^2$-factor only in a trivial way.
It remains to recall that a parallel spinor descends from $\R^2$ to a nontrivial 
quotient (flat cylinder or torus) if and only if the fundamental group acts on 
the spin structure of $\R^2$ in a trivial way, that is, the quotient $\Sigma^2$ 
carries the trivial spin structure.
\hfill$\square$

We end this section with the question -- asked by Ilka Agricola -- whether skew 
Killing spinors can be seen 
as parallel spinors w.r.t. a covariant derivative induced by some metric 
connection on $(TM,g)$.

\begin{prop}\label{p:sksparallel}
Let $(M^4,g)$ be any Riemannian spin manifold and $\psi$ be any nonzero skew 
Killing spinor on $M$.
Assume that, w.r.t. the splitting $\psi=\psi^++\psi^-$, both $\psi^\pm$ do not 
vanish on $M$.
Assume the existence of a metric connection $\nabla'$ on $(TM,g)$ such that 
$\psi$ is parallel w.r.t the covariant derivative induced by $\nabla'$ on 
$\Sigma M$.\\
Then $A\xi=0$, in particular 
$|\psi^+|=|\psi^-|$.
Moreover, 
$\nabla'_X=\nabla_X+2\left((AX\wedge\frac{\xi}{|\xi|^2}
)_+-(AX\wedge\xi)_-\right)$ for all $X\in TM$.
\end{prop}

{\it Proof}: Write $\nabla'=\nabla-B$ for some unknown $B\in 
T^*M\otimes\Lambda^2 T^*M$.
Recall that, for any $X\in TM$, 
$BX\in\mathrm{End}(TM)$ must be skew-symmetric because of both $\nabla,\nabla'$ 
being metric.
Then for any section $\varphi\in \Sigma M$ and any 
$X\in TM$,
\[\nabla_X'\varphi=\nabla_X\varphi-\frac{1}{2}BX\cdot\varphi,\]
where we see $BX$ as a two-form acting by Clifford multiplication on $\Sigma M$.
Since by assumption $\psi^+$ does not vanish anywhere, $\xi$ is a nowhere 
vanishing vector field on $M$.
The question is now whether $B$ exists such that 
\[\frac{1}{2}BX\cdot\psi= AX\cdot\psi\]
holds for all $X\in TM$.
Using the splitting $\psi=\psi^++\psi^-$, we obtain the following equivalent 
systems: 
\begin{eqnarray*}
\left\{\begin{array}{ll}\frac{1}{2}BX\cdot\psi^+&= AX\cdot\psi^-\\
&\\\frac{1}{2}BX\cdot\psi^-&= 
AX\cdot\psi^+\end{array}\right.&\Longleftrightarrow&\left\{\begin{array}{ll}
\frac{1}{2}BX\cdot\psi^+&= -AX\cdot\frac{\xi}{|\xi|^2}\cdot\psi^+\\
&\\\frac{1}{2}BX\cdot\psi^-&= AX\cdot\xi\cdot\psi^-\end{array}\right.\\
&\Longleftrightarrow&\left\{\begin{array}{ll}\frac{1}{2}BX\cdot\psi^+&= 
-(AX\wedge\frac{\xi}{|\xi|^2})\cdot\psi^++\langle 
AX,\frac{\xi}{|\xi|^2}\rangle\psi^+\\
&\\\frac{1}{2}BX\cdot\psi^-&= (AX\wedge\xi)\cdot\psi^--\langle 
AX,\xi\rangle\psi^-\end{array}\right.\\
&\Longleftrightarrow&\left\{\begin{array}{ll}\left(\frac{1}{2}BX+(AX\wedge\frac{
\xi}{|\xi|^2})\right)\cdot\psi^+&=\langle AX,\frac{\xi}{|\xi|^2}\rangle\psi^+\\
&\\\left(\frac{1}{2}BX-(AX\wedge\xi)\right)\cdot\psi^-&=-\langle 
AX,\xi\rangle\psi^-. \end{array}\right.
\end{eqnarray*}
Recall that a real $2$-form acts in a skew-Hermitian way on $\Sigma M$, 
therefore we obtain $\langle AX,\xi\rangle=0$ for all $X\in TM$ and thus 
$A\xi=0$.
Moreover, since self-dual resp. anti-self-dual $2$-forms kill negative resp. 
positive half spinors, the preceding systems gets equivalent to 
\[\left\{\begin{array}{ll}\left(\frac{1}{2}BX+(AX\wedge\frac{\xi}{|\xi|^2}
)\right)_+\cdot\psi^+&=0\\
&\\\left(\frac{1}{2}BX-(AX\wedge\xi)\right)_-\cdot\psi^-&=0 
\,.\end{array}\right.\]
On the other hand, as we have seen above, the maps
$
\bigwedge^2_-M \longrightarrow \Sigma^-M\cap (\psi^-)^\perp,\ 
\omega_-\longmapsto \omega_-\cdot \psi^-$ and  $\bigwedge^2_+M  \longrightarrow 
\Sigma^+M\cap (\psi^+)^\perp,\ \omega_+\longmapsto \omega_+\cdot \psi^+
$
are isomorphisms if $\psi^+\not=0$ and $\psi^-\not=0$.
Therefore we can deduce that $\left(\frac{1}{2}BX+(AX\wedge\frac{\xi}{|\xi|^2}
)\right)_+=0$ and $\left(\frac{1}{2}BX-(AX\wedge\xi)\right)_-=0$, which yields 
$BX=-2\left((AX\wedge\frac{\xi}{|\xi|^2}
)_+-(AX\wedge\xi)_-\right)$ and concludes the proof of Proposition 
\ref{p:sksparallel}.
\findemo

With other words, only a special subcase of the degenerate case can be 
considered with that ansatz, namely that considered in Proposition 
\ref{h=1/2}.
As a consequence, the general classification of 
$4$-dimensional Riemannian spin manifolds with skew Killing spinors cannot be 
obtained that way.
\color{black}

\section{Skew Killing spinors with non-degenerate Killing map~$A$}\label{S5}
This section is devoted to the case where we have a skew Killing spinor $\psi$ 
whose Killing map $A$ is non-degenerate everywhere. Recall that $\psi$ defines a 
vector field $\eta$ by (\ref{xieta}). As above, we put $\rho:=|\eta|$. Here, we 
want to assume that $M''=\left\{x\in M\mid 
\rho(x)\not\in\{0,1/2\}\right\}=\{x\in M\mid f(x)\not\in\{0,\pm1\} \}$ is equal 
to $M$. This is a sensible restriction since $M''$ is dense in $M$ if $A$ is 
non-degenerate everywhere, see Section~\ref{S3}. Working on $M''$ has the 
advantage that we do not have to care about the sign of $f$. Indeed, as 
explained in Remark~\ref{33}, up to a possible change of orientation on each 
connected component we may assume that $f>0$.  In particular, $f$ is defined by 
$\rho=|\eta|$ via $f=\sqrt{1-4\rho^2}$, which will be important for the reverse 
direction of Proposition~\ref{P1}.

\subsection{Equivalent description by complex structures}\label{S51}

Let $M$ be a manifold and $A$ be a skew-symmetric endomorphism field on $M$. 
Define a tensor field $C_A$ on $M$ by $C_A(X,Y):=(\nabla _XA)(Y)-(\nabla 
_YA)(X)$.

\begin{prop}\label{P1}
Let $M$ be a four-dimensional spin manifold and $A$ be a skew-symmetric 
endomorphism field on $M$. Put $C:=C_A$.

If $(M,g)$ admits a skew Killing spinor $\psi$ associated with $A$ such that 
$M=M''$, then there exist an almost Hermitian structure $J$ and a nowhere 
vanishing vector field $\eta$ of length $|\eta|=:\rho< 1/2$ such that 
\begin{eqnarray}
&&(\nabla_Y J)(X) = \frac{4}{f-1}X\lrcorner\, \big(J\eta\wedge AY +\eta\wedge 
JAY\big) \label{nJ},\\
&&\nabla \eta = f A \label{nxi},\\
&&g(C(\eta,X),J\eta)=\rho^2f\cdot g(C_P,X)
\label{CX}\\
&&g(C(J\eta,Z),J\eta)=*(C_P\wedge Z\wedge \eta\wedge J\eta),\quad Z\in P 
:=\{\eta,J\eta\}^\perp \label{star},
\end{eqnarray}
where $f:=\sqrt{1-4\rho^2}$ and $C_P:=C(s,Js)$ for any unit vector $s\in P$,
and such that the sectional curvature $K_P$ in direction $P$ satisfies
\begin{equation}
K_P= -\rho^{-2}g(C_P,J\eta)+4A_P^2, \label{s}
\end{equation}
where $A_P:=g(As,Js)$ for any unit vector $s\in P$.

If $M$ is simply-connected, then also the converse statement is true.
\end{prop} 

\begin{elemme} \label{L2} 
Assume that $J$, $A$ and $\eta$ satisfy Equations 
{\rm (\ref{nJ})} and {\rm (\ref{nxi})}. Then
\begin{eqnarray}
&&g(C(X,Y),\eta)=0,\label{gC}\\
&&R(X,Y)\eta=fC(X,Y)-4\eta\lrcorner (AX\wedge AY ),\label{Rxi}\\
&&R(X,Y)J\eta=-JC(X,Y)+\frac{4}{f-1}g(C(X,Y),J\eta)\eta-4J(\eta)\lrcorner 
(AX\wedge AY ) \label{RJxi}.
\end{eqnarray}
\end{elemme}
{\it Proof}: Note first that $X(f^2)=X(1-4|\eta|^2)=-8g(\nabla_X\eta,\eta) 
=-8fg(AX,\eta)$. This implies $X(f)=-4g(AX,\eta)$, which we will use in the 
following.
Let $X$ and $Y$ be vector fields on $M$ and assume that $\nabla X=\nabla Y=0$ 
holds at a point $p\in M$. At $p$, we have 
\begin{eqnarray*}
\nabla_X\nabla_Y \eta&=&\nabla_X \left(f AY\right) = -4g(AX,\eta)AY +f(\nabla_X 
A)Y.
\end{eqnarray*}
Thus 
$$R(X,Y)\eta=-4\big(g(AX,\eta)AY -g(AY,\eta)AX\big)+f C(X,Y),$$
which gives Eq.~(\ref{Rxi}).
In particular, this yields $0=R(X,Y,\eta,\eta)=g(C(X,Y),\eta)$, which proves 
(\ref{gC}) since $f\not=0$ everywhere.

In the following computation, the sign `$\equiv$' means equality up to a term 
$S(X,Y)$ for some symmetric bilinear map $S$. We compute
\begin{eqnarray*}
\lefteqn{\nabla_X\nabla_Y(J\eta)=\nabla_X\big((\nabla_Y 
J)(\eta)+J(\nabla_Y\eta)\big) }\\
&=& \nabla_X\Big(\frac4{f-1}\big( -g(\eta,AY)J\eta+\rho^2JAY- 
g(\eta,JAY)\eta\big)+J(\nabla_Y\eta)\Big)\\
&\equiv&\frac{16}{(f-1)^2}\,g(AX,\eta)\big( \rho^2JAY-g(\eta,JAY)\eta\big)\\
&&+\frac4{f-1}\Big(-g(\eta,(\nabla_XA)Y)J\eta-g(\eta,AY)(\nabla_XJ)\eta  
-fg(\eta,AY)JAX  +2 fg(AX,\eta)JAY \\[-1ex] 
&&\qquad \qquad \ +\rho^2(\nabla_XJ)(AY) +\rho^2J(\nabla_XA)Y  -fg(AX,JAY)\eta 
-g(\eta,(\nabla_XJ)AY)\eta\\
&&\qquad \qquad \  -g(\eta,J(\nabla_XA)Y)\eta -fg(\eta,JAY)AX \Big)\\
&&+f(\nabla_XJ)(AY)+J(\nabla_X\nabla_Y\eta)\\
&\equiv&\frac{16}{(f-1)^2}\,g(AX,\eta)\big( \rho^2JAY-g(\eta,JAY)\eta\big)\\
&&+\frac4{f-1}\Big(-g(\eta,(\nabla_XA)Y)J\eta-g(\eta,AY)(\nabla_XJ)\eta+2fg(AX,
\eta)JAY+\rho^2(\nabla_XJ)AY \\[-1ex]
&&\qquad \qquad \ +\rho^2J(\nabla_XA)Y -g(\eta,(\nabla_XJ)AY)\eta 
-g(\eta,J(\nabla_XA)Y)\eta -2fg(\eta,JAY)AX \Big)\\
&&+J(\nabla_X\nabla_Y\eta)\\
&\equiv&\frac{16\rho^2}{(f-1)^2}\,\big( g(AX,\eta)JAY+g(J\eta,AY)AX\big)\\
&&+\frac4{f-1}\Big(-g(\eta,(\nabla_XA)Y)J\eta+2fg(AX,\eta)JAY 
+\rho^2J(\nabla_XA)Y  \\[-1ex]
&&\qquad \qquad \ -g(\eta,J(\nabla_XA)Y)\eta -2fg(\eta,JAY)AX \Big)\\
&&+J(\nabla_X\nabla_Y\eta)\\
&=&4g(AX,\eta)JAY+4g(AY,J\eta)AX-\frac4{f-1}\big(g(\eta,(\nabla_XA)Y)J\eta 
+g(\eta,J(\nabla_XA)Y)\eta\big)\\
&&-(f+1)J(\nabla_XA)Y+J(\nabla_X\nabla_Y\eta).
\end{eqnarray*}
This implies 
\begin{eqnarray*}
\lefteqn{R(X,Y)J\eta\ =\ 4g(AX,\eta)JAY-4g(AY,\eta)JAX 
+4g(AY,J\eta)AX-4g(AX,J\eta)AY}\\
&&-\frac4{f-1}\big( g(\eta, C(X,Y))J\eta -g(J\eta,C(X,Y))\eta \big) 
-(f+1)JC(X,Y)+J(R(X,Y)\eta).
\end{eqnarray*}
Using Equations (\ref{gC}) and (\ref{Rxi}) we obtain (\ref{RJxi}).
\findemo

{\it Proof of Prop.\ref{P1}}: Before we start the proof of the two directions of 
the assertion, let us first suppose that, on $M$, we are given a Hermitian 
structure $J$ and a nowhere vanishing vector field $\eta$ of length $\rho<1/2$. 
We want to define a vector field $\xi$ such that the identities 
$\xi=-(|\xi|/\rho)\cdot\eta$ and $\rho=|\xi|/(1+|\xi|^2)$ hold according to 
Equation~(\ref{rhoxi}). Since this leads to a quadratic equation, we have to 
choose one of the solutions. Here we use our assumption $M=M''$
and define $f=\sqrt{1-4\rho^2}$ and
$\xi=2(f-1)^{-1}\eta$, compare Remark~\ref{33}, which motivates this choice. 
Assume that the orientation on $M$ is such that orthonormal bases of the form 
$s_1, Js_1,s_2,Js_2$ are negatively oriented. 
We define a one-dimensional subbundle $E$ of $\Sigma M$ by 
\begin{equation} E:= \{\ph \mid J(X)\cdot\ph^-=iX\cdot\ph^-,\ 
\ph^+=\xi\cdot\ph^-\}.\label{E} \end{equation}
We want to show that $E$ is parallel with respect to $\hat\nabla$ defined by 
$\hat \nabla_X \ph:=\nabla_X\ph-AX\cdot \ph$ if and only if $J$ and $\eta$ 
satisfy (\ref{nJ}) and (\ref{nxi}). Let $X$ and $Y$ be vector fields satisfying 
$\nabla X=\nabla Y=0$ at $p\in M$. Then we have at $p\in M$
\begin{eqnarray*}
\lefteqn{J(X)\cdot(\hat \nabla_Y\ph)^-=J(X)\cdot(\nabla_Y\ph -AY\cdot \ph)^-}\\
&=&J(X)\cdot(\nabla_Y\ph^- -AY\cdot \ph^+)\\
&=&\nabla_Y(J(X)\cdot\ph^-)-(\nabla_YJ)(X)\cdot\ph^-- J(X)A(Y)\cdot \ph^+\\
&=&\nabla_Y(iX\cdot\ph^-)-(\nabla_YJ)(X)\cdot\ph^-+ A(Y)J(X)\xi\cdot 
\ph^-+2g(JX,AY)\ph^+\\
&=&iX\cdot\nabla_Y\ph^--(\nabla_YJ)(X)\cdot\ph^-- iA(Y)\xi X\cdot 
\ph^--2g(JX,\xi)AY\cdot\ph^-+2g(JX,AY)\ph^+\\
&=&iX\cdot\nabla_Y\ph^--(\nabla_YJ)(X)\cdot\ph^-- iXA(Y)\xi \cdot \ph^- + 
2ig(\xi, X)AY\cdot\ph^- \\
&&- 2ig(AY,X)\xi\cdot\ph^--2g(JX,\xi)AY\cdot\ph^-+2g(JX,AY)\xi \cdot\ph^-\\
&=&iX\cdot\nabla_Y\ph^--(\nabla_YJ)(X)\cdot\ph^-- iXA(Y)\cdot\ph^+ + 2g(\xi, 
X)JA(Y)\cdot\ph^- \\
&&- 2g(AY,X)J(\xi)\cdot\ph^--2g(JX,\xi)AY\cdot\ph^-+2g(JX,AY)\xi \cdot\ph^-.
\end{eqnarray*}
This equals $iX\cdot(\hat\nabla_Y\ph)^-$ if and only if 
$(\nabla_Y J)(X) = 2X\lrcorner\, \big(J\xi\wedge AY +\xi\wedge JAY\big)$
holds, which is equivalent to Equation~(\ref{nJ}). Furthermore,
\begin{eqnarray*}
(\hat\nabla_X\ph)^+&=& \nabla_X\ph^+-AX\cdot\ph^-
\ =\ \nabla_X(\xi\cdot\ph^-)-AX\cdot\ph^-\\
&=&(\nabla_X\xi)\cdot\ph^-+\xi\cdot\nabla_X\ph^--AX\cdot\ph^-\\
&=&(\nabla_X\xi)\cdot\ph^-+\xi A(X)\cdot\ph^+  
+\xi\cdot(\hat\nabla_X\ph)^--AX\cdot\ph^- \\
&=&\left( 
\nabla_X\xi-(1-|\xi|^2)AX+2g(A\xi,
X)\xi\right)\cdot\ph^-+\xi\cdot(\hat\nabla_X\ph)^-.
\end{eqnarray*}
This equals $\xi \cdot(\hat\nabla_X\ph)^-$ if and only if $\nabla_X \xi = 
(1-|\xi|^2)AX-2g(A\xi,X)\xi$ holds, which is equivalent to (\ref{nxi}).  
Consequently,  $E$ is parallel with respect to $\hat\nabla$ if and only if $J$ 
and $\eta$ satisfy (\ref{nJ}) and (\ref{nxi}). 

Assume that $\hat \nabla$ reduces to a connection $\hat\nabla ^E$ on $E$. Then 
Equations (\ref{nJ}) and (\ref{nxi}), and therefore also (\ref{gC}), (\ref{Rxi}) 
and (\ref{RJxi}) hold. We will show that the curvature $\hat R$ of $\hat\nabla 
^E$ vanishes if and only if the Riemannian curvature $R$ of $M$ equals the 
tensor $B$ defined by
\begin{equation}\label{R}
B(X,Y):=\rho^{-2}\big( *(C(X,Y)\wedge \eta)- fC(X,Y)\wedge \eta\big) -4 AX\wedge 
AY
\end{equation}
for all vector fields $X$ and $Y$ on $M$.
By an easy calculation similar to that in the proof of 
Proposition~\ref{p:integrcond}, we get
\begin{eqnarray*}
\hat R_{X,Y}\ph&=&\textstyle{\frac12}R(X,Y)\cdot \ph -C(X,Y)\cdot\ph +2 
\big(AX\wedge AY\big)\cdot\ph.
\end{eqnarray*}
This shows that $\hat R$ vanishes if and only if 
\begin{equation} \label{hatR}
R(X,Y)\cdot \ph =2C(X,Y)\cdot\ph -4 \big(AX\wedge AY\big)\cdot\ph
\end{equation}
for all vector fields $X$ and $Y$ and all sections $\ph$ of $E$.
In the following, we will use that $\bigwedge_{\pm}^2M$ acts trivially on 
$\Sigma^{\mp} M$ and that, for any nowhere vanishing section $\ph^\pm$ of 
$\Sigma^\pm M$, the maps defined by (\ref{action}) are isomorphisms. Let $\ph$ 
be a section of $E$ such that $\ph^+(x)\not=0$, $\ph^-(x)\not=0$ for all $x\in 
M$ (here we use that $\xi$ does not vanish). 
Then 
\begin{eqnarray*}
2C(X,Y)\cdot\ph &=& 2C(X,Y)\cdot 
\big(\xi\cdot\ph^--|\xi|^{-2}\xi\cdot\ph^+\big)\\
&\bui{=}{(\ref{gC})} &2(C(X,Y)\wedge \xi)\cdot\ph^--2|\xi|^{-2}(C(X,Y)\wedge 
\xi)\cdot\ph^+\\
&=&2(C(X,Y)\wedge \xi)_-\cdot\ph-2|\xi|^{-2}(C(X,Y)\wedge \xi)_+\cdot\ph\\
&=&\frac4{f-1}(C(X,Y)\wedge \eta)_-\cdot\ph -\frac {f-1}{\rho^2} (C(X,Y)\wedge 
\eta)_+\cdot\ph\\
&=&\rho^{-2}\big(*(C(X,Y)\wedge \eta)- fC(X,Y)\wedge \eta \big)\cdot\ph.
\end{eqnarray*}
Thus (\ref{R}) and (\ref{hatR}) show that $\hat R$ vanishes if and only if 
$R=B$. The latter condition is equivalent to the system of equations
\begin{eqnarray}
R(X,Y)\eta&=&B(X,Y)\eta \label{Rxi2}\\
R(X,Y)J\eta&=&B(X,Y)J\eta\label{RJxi2}\\
R(s,Js,s,Js) &=& g (B(s,Js)s,Js)\label{Rs2}\\
g(B(\eta,X)s,Js)&=&g(B(s,Js)\eta,X)\label{B1}\\
g(B(J\eta,Z)s,Js)&=&g(B(s,Js)J\eta,Z)\label{B2}
\end{eqnarray}
for all $X,Y\in {\frak X}(M)$ and all $Z\in \Gamma(P)$. Recall that (\ref{gC}) 
holds in our situation, which we will use in the following computations.
Equations (\ref{B1}) and (\ref{B2}) are equivalent to the two equations
\begin{eqnarray*}
&g\big(s\lrcorner \big(*(C(\eta,X)\wedge \eta)-fC(\eta,X)\wedge \eta 
\big),Js\big)\ =\ g\big(\eta\lrcorner\big(*(C_P\wedge \eta)-fC_P\wedge 
\eta\big),X\big),&\\
&g\big(s\lrcorner \big(*(C(J\eta,Z)\wedge \eta)-fC(J\eta,Z)\wedge \eta 
\big),Js\big)\ =\ g\big(J\eta\lrcorner\big(*(C_P\wedge \eta)-fC_P\wedge 
\eta\big),Z\big),&
\end{eqnarray*}
which are equivalent to (\ref{CX}) and (\ref{star}), respectively. Because of
$$\eta\lrc\big( *(C(X,Y)\wedge \eta)- fC(X,Y)\wedge \eta\big)=f\rho^2 C(X,Y),$$
and
\begin{eqnarray*}
\lefteqn{\hspace{-1cm}J(\eta)\lrcorner \big( *(C(X,Y)\wedge \eta)- fC(X,Y)\wedge 
\eta\big) \ =\ *(C(X,Y)\wedge\eta\wedge J\eta) -f g(C(X,Y),J\eta)\eta}\\
&=& -\rho^2g(C(X,Y),s)Js+\rho^2g(C(X,Y),Js)s-f g(C(X,Y),J\eta)\eta\\
&=& -\rho^2\big(  g(JC(X,Y),Js)Js+g(JC(X,Y),s)s\big)-f g(C(X,Y),J\eta)\eta\\
&=&-\rho^2JC(X,Y)-(f+1)g(C(X,Y),J\eta)\eta\\
&=&-\rho^2JC(X,Y)+\frac{4\rho^2}{f-1}g(C(X,Y),J\eta)\eta,
\end{eqnarray*}
Lemma~\ref{L2} shows that Equation  (\ref{Rxi2}) is equivalent to (\ref{Rxi}) 
and (\ref{RJxi2}) is equivalent to (\ref{RJxi}). Recall that  (\ref{Rxi}) and 
(\ref{RJxi}) are satisfied in our situation.
Finally,
$$g \big(s\lrcorner \big(*(C_P\wedge \eta)- fC_P\wedge \eta\big),Js\big)
\ =\ g\big(*(s\wedge C_P\wedge\eta),Js\big) \ =\ g(C_P,J\eta),$$
which implies that (\ref{Rs2}) is equivalent to (\ref{s}).  Consequently, the 
curvature $\hat R$ of $\hat\nabla$ vanishes if and only if the Equations 
(\ref{CX}), (\ref{star}) and (\ref{s}) hold.

Now we can prove both directions of the proposition. Suppose that there exists a 
spinor field $\psi$ on $M$ satisfying $\nabla_X \psi = AX\cdot\psi$ for all 
$X\in TM$ such that $M=M''$. The latter condition means that the vector field 
$\eta$ defined in (\ref{xieta}) satisfies $0<\rho=|\eta|<1/2$. In particular, 
$\psi^-\not=0$ everywhere and we can define an almost Hermitian structure $J$ by 
$J(X)\cdot\psi^-=iX\cdot\psi^-$. Thus we may apply our above considerations. If 
we define $E\subset \Sigma M$ and $\hat\nabla$ as above, then $\psi$ is a 
$\hat\nabla$-parallel section of $E$. In particular, $\hat \nabla$ reduces to a 
connection $\hat \nabla^E$ and the curvature of $\hat \nabla^E$ vanishes thus 
(\ref{nJ}) -- (\ref{s}) hold.

Conversely, if we are given an almost Hermitian structure $J$ and a nowhere 
vanishing vector field $\eta$ of length $0<\rho=|\eta|<1/2$ such that (\ref{nJ}) 
-- (\ref{s}) are satisfied. Then we can define a one-dimensional subbundle 
$E\subset\Sigma M$ by (\ref{E}) together with a flat covariant derivative $\hat 
\nabla$ on $E$.  If $M$ is simply-connected, then $E$ admits a parallel section, 
which is a skew Killing spinor.
\findemo

\begin{re}\label{rint}
Let $J$ be an almost Hermitian structure on a four-dimensional manifold $M$ such 
that (\ref{nJ}) and (\ref{nxi}) hold for a skew-symmetric endomorphism field $A$ 
and a vector field $E$.  Then $J$ defines a reduction of the  $\SO(4)$-bundle 
$\SO(M)$ to $\U(2)$. Here we want to give the intrinsic torsion of this bundle 
in the special case where $A$ and $J$ commute. The two components of the 
intrinsic torsion of this bundle are the Nijenhuis tensor $N$ of $J$ and the 
differential $d\Omega$ of the K\"ahler form $\Omega:=g(J\cdot, \cdot)$. A direct 
calculation using (\ref{nJ}) and (\ref{nxi})  shows that under the assumption 
$AJ=JA$ these components are given by $N=0$ and  $d\Omega=-2A\wedge(\xi\lrcorner 
\Omega)$.
\end{re}

\subsection{The case where $A\eta$ is parallel to $J\eta$}\label{S52}
Let us assume again that the Killing map $A$ is non-degenerate everywhere.  We 
want to consider the case where $A\eta$ is parallel to $J\eta$ in more detail. 
We will see that, in this situation, the existence of skew Killing spinors is 
related to doubly warped products and to local DWP-structures. These notions and 
their basic properties are explained in the appendix. 

\begin{elemme} Assume that $M$ admits a skew Killing spinor with nowhere 
vanishing Killing map $A$ that satisfies  $A\eta=u J\eta$ for some function $u.$ 
Then $A^2\eta=-u^2\eta.$ In particular, $AJ=JA.$ 
\end{elemme}
{\it Proof.} Note first that Lemma~\ref{l:h}, \ref{i:eta4} and Eq.~(\ref{gC}) 
give
\begin{eqnarray*}
0&=&f(dA)(X,Y,\eta)\ =\ f\big( 
g((\nabla_XA)Y,\eta)-g((\nabla_YA)X,\eta)+g((\nabla_ \eta A)X,Y)\big)\\
&=&f g(C(X,Y),\eta)+fg((\nabla_\eta A)X,Y)\ =\ fg((\nabla_\eta A)X,Y)
\end{eqnarray*}
for all $X,Y\in TM$. Consequently, $f\nabla_\eta A=0$.
Because of  
$$\eta\lrc (J\eta\wedge A\eta +\eta\wedge 
JA\eta)=|\eta|^2JA\eta-g(\eta,JA\eta)\eta=-u|\eta|^2\eta +u|\eta|^2\eta=0,$$ 
Eq.~(\ref{nJ}) gives $(\nabla_\eta J)\eta=0$. Now, by differentiating the 
equality $A\eta=u J\eta$ in the direction of $\eta$, we get 
$$
\nabla_\eta A\eta=(\nabla_\eta A)\eta+A(\nabla_\eta\eta)=\eta(u)J\eta+u 
(\nabla_\eta J)\eta+u J(\nabla_\eta\eta) =
\eta(u)J\eta+u J(\nabla_\eta\eta).
$$
Finally, using the fact that $\nabla_\eta\eta=fA\eta$ and $f\nabla_\eta A=0,$ we 
get that $\eta(u)=0$ and $f^2A^2\eta=-u^2f^2\eta.$ The latter equation implies  
$A^2\eta=-u^2\eta$ since supp$(f)=M$. 
\findemo

Let $(\hat M^3,\hat g, \hat \eta)$ be a minimal Riemannian flow, i.e., an 
orientable three-dimensional Riemannian manifold together with a unit Killing 
vector field $\hat \eta$. Then, locally, $(\hat M,\hat g)$ is a Riemannian 
submersion over a two-dimensional base manifold $B$. Let us fix a Hermitian 
structure $\hat J$ on $\hat\eta^\perp$ and put $\omega:=\hat g(\cdot,\hat 
J\cdot)$. We define a function $\hat\tau$ on $\hat M$ which is constant along 
the fibres by $\hat \nabla_X\hat\eta=\hat\tau\cdot \hat J(X)$ for $X\in\hat 
\eta^\perp$. Furthermore, let $\hat K$ denote the Gaussian curvature of $B$. Now 
consider the metric $g_{rs}=r^2\hat g_{\hat \eta} \oplus s^2\hat 
g_{\hat\eta^\perp}$ on $\hat M$,  where $\hat{g}_{\hat{\eta}}$, 
$\hat{g}_{\hat{\eta}^\perp}$ are the components of the metric $\hat{g}$ along 
$\R\hat{\eta}$ and $\hat{\eta}^\perp$, respectively.  Then $(\hat 
M,g_{rs},r^{-1}\hat\eta)$ is again a minimal Riemannian flow and we obtain new 
functions  $\hat \tau$ and $\hat K$, say $\hat \tau_{rs}$ and $\hat K_{rs}$. 
These functions satisfy
\begin{equation}\label{taurs}
\hat \tau_{rs} = rs^{-2}\hat \tau,\quad \hat K_{rs}=s^{-2}\hat K.
\end{equation}
If our four-dimensional manifold $M$ is endowed with a DWP-structure, then every 
three-dimensional leaf associated with this structure can be understood as a 
minimal Riemannian flow. In this way, we obtain functions $\tau$ and $K$ on 
$M$. 

\begin{ethm}\label{thm}
Assume that $M$ admits a skew Killing spinor such that $A\eta || J\eta$ and 
that 
$\rho=|\eta|\not\in\{0,1/2\}$ everywhere. Then $(\nu:=-\rho^{-1} J\eta,\eta)$ 
is 
a local DWP-structure on $M$ such that
\begin{equation}
f\cdot \mu=\tau,
\quad K=2\mu\lambda +2\tau ^2, \label{Ktau}
\end{equation}
 for $f:=\sqrt{1-4\rho^2}$, where $\lambda$ and $\mu$ are the eigenvalues of the 
Weingarten map $W=-\nabla\nu$ on $\R{\eta}$ and ${\eta}^\perp\cap\nu^\perp$, 
respectively. 

Conversely, suppose that $M$ is simply-connected and admits a local 
DWP-structure $(\nu,\eta)$ on $M$ such that the length $\rho$ of $\eta$ 
satisfies $0<\rho<1/2$.  Moreover, assume that $K$ and $\tau$ satisfy 
(\ref{Ktau}) for  $f:=\sqrt{1-4\rho^2}$. Then $M$ admits a skew Killing spinor 
such that $\eta$ is associated with $\psi$ according to (\ref{xieta}) and such 
that $A\eta || J\eta$.
\end{ethm}
{\it Proof}: Assume first that $M$ admits a skew Killing spinor such that 
$A\eta 
|| J\eta$ and $0<\rho<1/2$ everywhere. We define a vector field $\nu$ and 
functions $A_E$ and $A_P$ by
$$\nu=-\rho^{-1} J\eta,\quad AJ\eta=-A_E\eta,\quad AJZ=-A_PZ,\ Z\in 
\{\eta,\nu\}^\perp. $$
Then $\eta$ is a Killing vector field, see Remark~\ref{rxi}. Equation 
(\ref{nxi}) yields
\begin{equation}
\nu(\rho)= f A_E. \label{dgl1}
\end{equation}
We want to show that $(\nu,\eta)$ is a DWP-structure. The next Lemma will prove 
all properties of such a structure except the conditions for the Weingarten map 
$W=-\nabla\nu$ and its eigenvalues. 
\begin{elemme} Assume that $M$ admits a skew Killing spinor such that $A\eta || 
J\eta$ and $|\eta|\not\in\{0,1/2\}$  everywhere. Then
\begin{enumerate}
\item $\nu^\perp$ is integrable, 
\item the vector field $\eta$ has constant length on the integral manifolds of 
$\nu^\perp$,
\item the unit vector field $\nu$ is geodesic.
\end{enumerate}
\end{elemme}
{\it Proof}:  Take $X,Y \perp J\eta$. Using $JA=AJ$ we obtain
\begin{eqnarray*}
\lefteqn{g([X,Y],J\eta) \ = \ g(\nabla_XY,J\eta)-g(\nabla_YX,J\eta)\ =\ 
-g(Y,\nabla_X(J\eta))+g(X,\nabla_Y(J\eta))}\\
&=&-g(Y,(\nabla_XJ)\eta)-g(Y,J(\nabla_X\eta))  
+g(X,(\nabla_YJ)\eta)+g(X,J(\nabla_Y\eta))\\
&=&-4(f-1)^{-1} g\big(Y,\eta\lrcorner\, \big(J\eta\wedge AX +\eta\wedge 
JAX\big)\big)   -fg(Y, JAX)\\
&&+4(f-1)^{-1} g\big(X,\eta\lrcorner\, \big(J\eta\wedge AY +\eta\wedge 
JAY\big)\big)    + fg(X,JAY)\\
&=& 4(f-1)^{-1}\left(-g\big(Y,\rho^2 JAX-g(\eta,JAX)\eta\big) +g\big(X,\rho^2 
JA(Y)-g(\eta,JAY)\eta\big)\right) \\
&=&0
\end{eqnarray*}
since $JA\eta$ is a multiple of $\eta$. This proves the first claim. 
For $X\perp J\eta$, we have
\begin{eqnarray*}
Xg( \eta,\eta) = 2g(\nabla _X\eta,\eta) = 2fg(AX,\eta) =0
\end{eqnarray*}
since $A\eta || J\eta$. This shows the second assertion. The third one follows 
from (\ref{nJ}), (\ref{nxi}) and (\ref{dgl1}).
\findemo

We compute the eigenvalues of the Weingarten map $-\nabla\nu$, where we use that 
$\rho=|\eta|$ is constant on the integral manifolds of $\nu^\perp$: 
\begin{eqnarray}
-\nabla_\eta\nu &=&\rho^{-1} \nabla_\eta(J\eta)\ =\ \rho^{-1}(\nabla_\eta 
J)(\eta)+\rho^{-1} J(\nabla_\eta\eta)\nonumber\\
&=&4(f-1)^{-1}\rho^{-1} \cdot \eta\lrcorner (J\eta\wedge A\eta+\eta\wedge 
JA\eta)+f\rho^{-1}JA\eta\nonumber\\
&=&-f\rho^{-1}A_E\cdot\eta ,\label{nu1}
\end{eqnarray}
\begin{eqnarray}
-\nabla_Z\nu &=&\rho^{-1} \nabla_Z(J\eta)\ =\ \rho^{-1}(\nabla_Z 
J)(\eta)+\rho^{-1} J(\nabla_Z\eta) \nonumber\\
&=&4(f-1)^{-1}\rho^{-1} \cdot\eta \lrcorner (J\eta\wedge AZ+\eta\wedge 
JAZ)+f\rho^{-1}JAZ\nonumber \\
&=&\rho^{-1}A_P\cdot Z \label{nu2}
\end{eqnarray}
for $Z\in\nu^\perp\cap \eta^\perp$. Thus 
\begin{equation}\label{lambdamu}
\lambda=-f\rho^{-1}A_E,\quad \mu=\rho^{-1}A_P
\end{equation}
are the eigenvalues of $-\nabla\nu$. 
We fix a local section $s$ in $P=\{\eta,\nu\}^\perp$  and put
\begin{equation} s_1:=-\eta/\rho,\quad s_2:=Js_1=\nu,\quad s_3:=s,\quad s_4:=Js. 
\label{frame} \end{equation}
Let $s^1,\dots, s^4$ denote the dual local basis of $T^*M$.
By (\ref{nxi}), (\ref{nu1}) and (\ref{nu2}), the coefficients $\theta_{ij}:= 
g(\nabla s_i,s_j)$ of the Levi Civita connection satisfy
\begin{equation}\label{conncoeff}
\begin{array}{lll}\theta_{12}= -f\rho^{-1} A_E\, s^1,\quad & 
\theta_{13}=f\rho^{-1}A_P\,s^4,\quad &\theta_{14}=-f\rho^{-1}A_P\,s^3 , \\[1ex]
\theta_{23}= -\rho^{-1}A_P\,s^3, \quad &
\theta_{24}= -\rho^{-1}A_P\,s^4\,.&
\end{array}
\end{equation}
This gives
$$C_P=-2\rho^{-1}\big( A_P^2 + f A_P A_E \big) s_2 - s_3(A_P)\cdot s_3-s_4(A_P) 
\cdot s_4.$$
Indeed, (\ref{gC}) shows that $g(C_P,s_1)=0$. Furthermore,
\begin{eqnarray*}
g(C_P,s_2)&=& g((\nabla_{s_3}A)(s_4)-(\nabla_{s_4}A)(s_3), s_2)\\
&=& s_3\big(g(As_4,s_2)\big)-g(A(\nabla_{s_3}s_4),s_2)-g(As_4,\nabla_{s_3}s_2)\\
&&-s_4\big(g(As_3,s_2)\big)+g(A(\nabla_{s_4}s_3),s_2)+g(As_3,\nabla_{s_4}s_2)\\
&=&g(\nabla_{s_3}s_4,As_2)-g(As_4,\nabla_{s_3}s_2)-g(\nabla_{s_4}s_3,
As_2)+g(As_3,\nabla_{s_4}s_2)\\
&=&A_E\big(\theta_{14}(s_3)-\theta_{13}(s_4)\big)+A_P\big(\theta_{23}
(s_3)+\theta_{24}(s_4)\big),\\[2ex]
g(C_P,s_3)&=&g((\nabla_{s_3}A)(s_4)-(\nabla_{s_4}A)(s_3), s_3)\ =\ 
g((\nabla_{s_3}A)(s_4), s_3)\\
&=&g(\nabla_{s_3}(As_4)-A(\nabla_{s_3}s_4), s_3)\\
&=&-g(\nabla_{s_3}(A_Ps_3),s_3)+g(\nabla_{s_3}s_4,A_Ps_4)\\
&=&-s_3(A_P).
\end{eqnarray*}
Analogously, $g(C_P,s_4)=-s_4(A_P)$. Equations~(\ref{conncoeff}) imply
$$g(C(J\eta,Z),J\eta)=\rho^2g(C(s_2,Z),s_2)=\rho^2\big(s_2(g(AZ,s_2))- 
g(A(\nabla_{s_2}Z),s_2)-g(AZ,\nabla_{s_2}s_2)\big)=0$$
for $Z\in\{s_3,s_4\}$ and, similarly,
\begin{eqnarray*}
g(C(\eta,X),J\eta)&=&\rho^2 \big( s_1(A_E g(X,s_1)) - A_E 
g(\nabla_{s_1}X,s_1)-g(AX,\nabla_{s_1}s_2)-X(A_E)\big)\\
&=&\rho^2\big(s_1(A_E)g(X,s_1)+A_E 
g(X,\nabla_{s_1}s_1)+g(X,A(\nabla_{s_1}s_2))-X(A_E)\big)\\
&=&\rho^2\big(s_1(A_E)g(X,s_1)-\frac{fA_E^2}{\rho}g(X,s_2)+\frac{fA_E^2}{\rho}
g(X,s_2)-X(A_E)\big)\\
&=&-\rho^2X(A_E)
\end{eqnarray*}
for $X\in\{s_2,s_3,s_4\}$. Furthermore, $g(C(\eta,s_1),J\eta)=0$ since $C$ is 
antisymmetric. 
Hence, under the assumption that (\ref{conncoeff}) holds, Eqs.~(\ref{CX}), 
(\ref{star}) and (\ref{s}) are equivalent to the system of equations
\begin{eqnarray}
&&s_2(A_E)= 2f\rho^{-1}A_P^2+2f^2\rho^{-1} A_E A_P \label{dgl2}\\
&&s_j(A_P)= s_j(A_E) \ =\ 0,\ \ j=3,4, \label{dgl}\\
\label{KP}
&&K_P=-2 f\rho^{-2}A_E A_P-2(\rho^{-2}-2)\,A_P^2 .
\end{eqnarray}
We also have
$$s_1(A_E)=s_1(A_P)=0.$$
Indeed, (\ref{gC}) implies $g\big(C(s_1,s_2),\eta\big)=0$, thus we obtain
$$0=g\big((\nabla_{s_1}A)s_2 -(\nabla_{s_2}A)s_1,\eta\big) 
=g\big((\nabla_{s_1}A)s_2,\eta\big)=g\big(\nabla_{s_1}(As_2),\eta\big)=\rho 
s_1(A_E),$$
which gives $s_1(A_E)=0$. Using (\ref{dgl2}) and taking into account that 
$[s_1,s_2]$ is a multiple of $s_1$, we get
$$
0=s_1(s_2(A_E))=2 s_1\big(f\rho^{-1}A_P^2+f^2\rho^{-1} A_E A_P\big) 
=2f\rho^{-1}\big(2A_P+f A_E \big) s_1(A_P).
$$
Assume that $s_1(A_P)(x)\not=0$ at $x\in M$. Then $s_1(A_P)\not=0$ in an open 
neighbourhood $U$ of $x$. But then $2A_P=-f A_E$ on $U$, which would imply 
$s_1(A_P)=0$, a contradiction.

Hence we proved that besides $\rho$ also $A_E$ and $A_P$ are constant on the 
integral manifolds of $\nu^\perp$. Thus also $\mu$ and $\lambda$ are constant 
along these leaves. Consequently, $(\nu, \eta)$ is a local DWP-structure on $M$. 
By (\ref{lambdamu}), the associated function $\tau$ satisfies 
$$
\tau=\rho^{-1}g(\nabla_s\eta,Js)
=\rho^{-1}g(fAs,Js)=f\mu,$$
where $s\in \{\eta,\nu\}^\perp$ is of length one.
This proves  the first equation in (\ref{Ktau}). 

It remains to prove that also the second equation in (\ref{Ktau}) is true. Let 
$N$ be an integral manifold of $\nu^\perp$. Then, locally, $N$ is a Riemannian 
submersion over a base manifold $B$. The following lemma will relate the 
sectional curvature $K_P$ in direction of $P=\Span\{s_3,s_4\}$ to the Gaussian 
curvature $K$ of $B$, which will almost finish the proof of the forward 
direction of  Theorem~\ref{thm}.
 
\begin{elemme}\label{LK} Let $(\nu, \eta)$ be a local DWP-structure such that 
the coefficients of the Levi-Civita connection satisfy (\ref{conncoeff}) with 
respect to an orthonormal frame $s_1=-\eta/\rho$, $s_2=\nu$, $s_3, s_4$. Then 
the Gaussian curvature $K$ of $B$ equals $$K=K_P+(1+3f^2)\rho^{-2}A_P^2.$$ 
\end{elemme}
{\it Proof}:  The second fundamental form $\alpha$ of $N\subset M$ satisfies
$$
\alpha(s_3,s_3)=\alpha(s_4,s_4) =\rho^{-1}A_P s_2, \quad
\alpha(s_3,s_4)=0, 
$$
which follows from (\ref{conncoeff}).
Hence the Gauss equation gives
\begin{eqnarray}
K_P&=&R(s_3,s_4,s_4,s_3)\ =\ 
R^{N}(s_3,s_4,s_4,s_3)-g\big(\alpha(s_3,s_3),\alpha(s_4,s_4)\big)\nonumber \\
&=& R^{N}(s_3,s_4,s_4,s_3)-\rho^{-2}A_P^2. \label{hilf}
\end{eqnarray}
Let ${\cal A} $ denote the fundamental tensor used in O'Neill's formulas. We 
have
$${\cal A}_{s_j}s_j = g(\nabla _{s_j}s_j, s_1)s_1 =0, \quad j=3,4$$
and
$${\cal A}_{s_3}s_4 =-{\cal A}_{s_4}s_3= g(\nabla _{s_3}s_4, s_1)s_1 = 
-\theta_{14}(s_3)s_1= f\rho^{-1}A_P s_1.$$
The O'Neill formula for $R^N$ now gives 
$$R^{N}(s_3,s_4,s_4,s_3)=K-3|{\cal A}_{s_3}s_4|^2
=K-3f^2\rho^{-2}A_P^2,
$$
which combined with (\ref{hilf}) implies the assertion. \findemo

Lemma \ref{LK} together with (\ref{lambdamu}) and (\ref{KP}) finally shows that 
$B$ has constant curvature
$$ K=- 2f\rho^{-2}A_EA_P+2f^2\rho^{-2}A_P^2= 2\mu\lambda +2\tau ^2.  $$

Now let $M$ be simply-connected and let $(\nu,\eta)$ be a local DWP-structure on 
$M$ such that $0<|\eta|<1/2$ and such that (\ref{Ktau}) holds. Note that 
$f=\sqrt{1-4\rho^2}$ is smooth since $\rho=|\eta|<1/2$.  By assumption, $\rho$ 
is constant on the integral leaves of $\nu^\perp$. We write $\partial_t$ for the 
derivative in direction $\nu$. By (\ref{dwp}), we have $\rho'=-\lambda \rho$, 
which implies
\begin{equation}\label{fprime} f'= 4\lambda\rho^2/f = 
\lambda(-f+1/f).\end{equation}

We define functions
\begin{equation}\label{func}
A_E:=-\lambda \rho f^{-1},\quad A_P:=\mu\rho,
\end{equation}
which are all constant along the integral leaves of $\nu^\perp$.
We consider a local orthonormal frame 
$$s_1:=-\eta/\rho,\quad s_2=\nu, \quad s_3,\quad s_4$$
such that $s_3,s_4$ is a positively oriented basis of $\{\eta,\nu\}^\perp$. 
The assumption that $(\nu,\eta)$ is a local DWP-structure with eigenvalues 
$\lambda$ and $\mu$ together with the assumption 
$\tau=f\mu$
 implies that the local coefficients of the Levi-Civita connection satisfy 
Equations~(\ref{conncoeff}). Indeed, 
$$\nabla_{s_2}s_2=0,\quad \nabla_{s_1}s_2=-\lambda s_1,\quad 
\nabla_{s_j}s_2=-\mu s_j, \ j=3,4$$
implies 
\begin{eqnarray*}
\theta_{12}&=&-g(\nabla {s_2},s_1)=-g(\nabla_{s_1}s_2,s_1)s^1=\lambda s^1,\\
\theta_{23}&=&g(\nabla {s_2},s_3)=g(\nabla_{s_3}s_2,s_3)s^3=-\mu s^3,
\end{eqnarray*}
and (\ref{func}) gives the formulas for $\theta_{12}$ and $\theta_{23}$. 
Similarly, we get $\theta_{24}$. On $\Span\{s_3,s_4\}$, we fix the Hermitian 
structure $J$ that maps $s_3$ to $s_4$. Recall that $\tau$ is defined by 
$\nabla_X\eta=\tau JX$ for all $X\in\Span\{s_3,s_4\}$.  Since $\eta$ is a 
Killing vector field and $\rho$ is constant along the integrals leaves, we 
obtain
\begin{eqnarray*}
\theta_{13}&=& g(\nabla_{s_1}s_1,s_3)s^1+ \dots +g(\nabla_{s_4}s_1,s_3)s^4\\
&=&-\rho^{-1} \big(\,g(\nabla_{s_1}\eta,s_3)s^1+g(\nabla_{s_2}\eta,s_3)s^2  
+g(\nabla_{s_4}\eta,s_3)s^4\big) \\
&=&-\rho^{-1} \big(-g(s_1,\nabla_{s_3}\eta)s^1-g(s_2,\nabla_{s_3}\eta)s^2  
+g(\nabla_{s_4}\eta,s_3)s^4\big) \\
&=&
 -g(s_1,\nabla_{s_3}s_1)s^1-g(s_2,\nabla_{s_3}s_1)s^2 +\tau s^4
\\
 &=& f\mu s^4\ =\ f\rho^{-1}A_P s^4,
\end{eqnarray*}
where we used the already proven equation $\theta_{12}(s_3)=0$.
 Analogously, we obtain $\theta_{14}$.
Now we define skew-symmetric maps $A$ and $J$ by 
\begin{eqnarray*}
&A(s_1)=A_Es_2,\quad A(s_2)=-A_Es_1,\quad A(s_3)=A_Ps_4,\quad  A(s_4)=-A_Ps_3,  
&\\
&J(s_1)=s_2,\quad J(s_2)=-s_1,\quad J(s_3)=s_4, \quad J(s_4)=-s_3.&
\end{eqnarray*}
Note that $J$ extends the above defined map $J$ on $\Span\{s_3,s_4\}$. A few 
lines above, we proved that (\ref{conncoeff}) holds in our situation. Using this 
equation, we obtain 
\begin{eqnarray*}
\nabla_{s_1}\eta&=&-\rho\nabla_{s_1}s_1\ =\ fA_Es_2\ =\ fA(s_1),\\
\nabla_{s_2}\eta&=&-\rho' s_1-\rho\nabla_{s_2}s_1\ =\ \lambda\rho s_1\ =\ -fA_E 
s_1 \ =\ fA(s_2),\\
\nabla_{s_3}\eta&=&-\rho\nabla_{s_3}s_1\ =\ fA_Ps_4\ =\ fA(s_3),\\
\nabla_{s_4}\eta&=&-\rho\nabla_{s_4}s_1\ =\ -fA_Ps_3\ =\ fA(s_4).
\end{eqnarray*}
Hence, $\eta$ satisfies (\ref{nxi}). By definition of $J$ and $A$, 
Eq.~(\ref{nJ}) is equivalent to the system of equations
$$\begin{array}{ll}
\nabla_\eta J=\nabla_{J\eta}J=0,&\\[1ex]
(\nabla_sJ)(\eta)=(f+1)A_Ps,\quad &(\nabla_sJ)(J\eta)=-(f+1)A_PJ(s),\\[1ex]
(\nabla_sJ)(s)=4(f-1)^{-1}A_P\eta,\quad  
&(\nabla_sJ)(Js)=-4(f-1)^{-1}A_PJ(\eta),
\end{array}$$
for all $s\in\{\eta,J\eta\}^\perp$, $|s|=1$, which indeed can be verified 
using~(\ref{conncoeff}). Finally, we prove that (\ref{CX}), (\ref{star}) and 
(\ref{s}) hold. We already have seen that these equations are equivalent to 
(\ref{dgl2}), (\ref{dgl}) and (\ref{KP}). Now we use Lemma~\ref{LK}. Together 
with our assumption (\ref{Ktau}) and Equation~(\ref{func}), it implies 
$$K_P=2\mu\lambda+2\tau^2-\rho^{-2}(1+3f^2)A_P^2=-2f\rho^{-2}A_PA_E-\rho^{-2}
(1+f^2)A_P^2 ,$$
which is equivalent to (\ref{KP}). Also (\ref{dgl}) holds since $\rho$, 
$\lambda$ and $\mu$ are constant on the leaves by assumption. It remains to 
prove~(\ref{dgl2}).
Locally, $(M,g)$ is isometric to a doubly warped product $(I\times\hat M, 
\rho(t)^2\hat g_{\hat \eta} \oplus \sigma(t)^2\hat g_{\hat\eta^\perp})$. In 
particular, $\sigma'= -\mu\sigma$ by~(\ref{dwp}). Furthermore, $\tau=\hat \tau 
\rho\sigma^{-2}$ and $K=\hat K \sigma^{-2}$ by (\ref{taurs}) for some constants 
$\hat \tau$ and $\hat K$. Thus, by assumption (\ref{Ktau}),
$$\mu f=\rho\sigma^{-2}\hat\tau.$$
Taking the absolute value and then the logarithm on both sides and 
differentiating, we obtain
$$\frac{\mu'}{\mu}+\frac 
{f'}f=\frac{\rho'}\rho-2\frac{\sigma'}\sigma=-\lambda+2\mu$$
and therefore, by (\ref{fprime}),
$$\mu'=(1-f^{-2})\lambda\mu -\lambda\mu+2\mu^2=-f^{-2}\lambda\mu+2\mu^2$$
holds (globally) on $M$.
By assumption, 
$$\hat K=K\sigma^2=(2\mu\lambda+2\tau^2)\sigma^2=2\mu(\lambda+\mu 
f^2)\sigma^2.$$
Differentiating, using $\sigma'=-\mu\sigma$ and dividing by $2\mu\sigma^2$ 
yields
\begin{eqnarray*}0&=&\frac{\mu'}\mu(\lambda+\mu f^2)+\lambda'+\mu'f^2+2\mu 
ff'-2\mu\lambda-2\mu^2 f^2\\
&=&(-\lambda f^{-2}+2\mu)(\lambda+\mu f^2)+\lambda'+(-\lambda\mu 
f^{-2}+2\mu^2)f^2+2(1-f^2) \lambda\mu-2\mu\lambda-2\mu^2 f^2,
\end{eqnarray*}
thus
$\lambda'=\lambda^2 f^{-2}-2f^2(\mu^2-\lambda\mu),$
which gives~(\ref{dgl2}) by Equations (\ref{fprime}) and (\ref{func}). 
Consequently, we proved that Equations (\ref{nJ}) -- (\ref{s}) hold. Now 
Proposition~\ref{P1} shows the existence of a skew Killing spinor.
\findemo
\begin{ecor}\label{corfinal}
Let $(M,g)$ admit a skew Killing spinor such that  $A\eta || J\eta$ and 
$|\eta|\not\in\{0,1/2\}$ everywhere. Then $M$ is locally isometric to a doubly 
warped product 
$(I\times \hat 
M,dt^2\oplus\rho(t)^2\hat{g}_{\hat{\eta}}\oplus\sigma(t)^2\hat{g}_{\hat{\eta}
^\perp})$ for which the data $\hat K$ and $\hat \tau$ are constant and $\rho$ 
and $\sigma$ satisfy the differential equations
\begin{eqnarray}
\label{sol1}
&&(\sigma^2)'=-\frac 2{\sqrt{1-4\rho^2}}\rho\hat \tau
\\
\label{sol2}
&&(\sigma^2)'\,\frac{\rho'}\rho=\hat K-2\frac{\rho^2}{\sigma^2}\hat\tau^2.
\end{eqnarray}

Conversely, if $M$ is isometric to a doubly warped product 
$(I\times \hat 
M,dt^2\oplus\rho(t)^2\hat{g}_{\hat{\eta}}\oplus\sigma(t)^2\hat{g}_{\hat{\eta}
^\perp})$  for which the data $\hat K$ and $\hat \tau$ are constant and $\rho$ 
and $\sigma$ satisfy the differential equations (\ref{sol1}) and (\ref{sol2}) 
and if  $\hat M$ is simply-connected, then $(M,g)$ admits a skew Killing spinor 
such that $A\eta || J\eta$.
\end{ecor}
{\it Proof}: 
The condition 
$\mu\cdot f=\tau$ is equivalent to $-\frac {\sigma'}\sigma\cdot 
f=\frac\rho{\sigma^2}\hat\tau$,
 thus to (\ref{sol1}),
and $K=2\mu\lambda +2\tau ^2$ is equivalent to $\frac{\hat 
K}{\sigma^2}=2\frac{\rho'}\rho\frac{\sigma'}\sigma+2\left(\frac\rho{\sigma^2}
\hat\tau\right)^2$, thus to (\ref{sol2}).
\findemo

\begin{erem}
Locally, Equations (\ref{sol1}) and (\ref{sol2}) can be solved explicitly to get 
solutions $\sigma$ and $\rho$.
\end{erem}
\begin{erem}\label{notequasi}
Let us study the restriction of a skew Killing spinor $\psi$ to $N$. The 
restriction $(\Sigma M)|_N$ can be understood using an isomorphism
$$\phi:\ (\Sigma M)|_N\longrightarrow \Sigma N \oplus\Sigma N=\phi((\Sigma^+ 
M)|_N)\oplus \phi((\Sigma^- M)|_N)$$
which is compatible with the Clifford multiplication in the following sense. If 
$\phi(\ph)=(u,v)$, then 
$$\phi(\nu\cdot\ph)=(-v,u),\quad \phi(\nu\cdot 
X\cdot\ph)=(-X\cdot_{\scriptscriptstyle N} u,X\clN v),$$
where $\nu=s_2$ is a normal vector of $N$, $X$ is a tangent vector of $N$ and 
`$\,\clN$' denotes the Clifford multiplication on $\Sigma N$. In particular, 
$s_1s_3s_4\clN u=u$ for all $u\in\Sigma N$. By the spinorial O'Neill formulas, 
we obtain 
\begin{equation}\label{Eres}
\nabla^N_{\eta}\phi(\psi^\pm)=-\frac\lambda{2f}\eta\clN\phi(\psi^{\pm}),\quad 
\nabla^N_{Z}\phi(\psi^\pm)= -\frac{\mu f}2Z\clN\phi(\psi^\pm)
\end{equation}
for all $Z\in TN\cap\eta^\perp$. Up to rescaling, these are Sasakian 
quasi-Killing spinors on $N$, which we will explain in the following.

Up to rescaling of the metric, each integral manifold $N$ in our construction 
has a Sasakian structure, see \cite{Blair} for a definition of such structures. 
Indeed, $\eta$ restricted to $N$ is a Killing vector field of constant length 
and $\nabla\eta$ restricted to $\eta^\perp$ equals $|\eta|\tau J|_{\eta^\perp}$, 
where also $\tau$ is constant. The Nijenhuis tensor of $J|_{\eta^\perp}$ 
vanishes since $\eta^\perp$ is two-dimensional. Consequently, $\tilde \xi:= 
\eta/(\tau |\eta|)$ is the Reeb vector field of a Sasakian structure on 
$(N,\tilde g:=\tau^2g)$.  The scalar curvature of $(N,\tilde g)$ equals 
$\S=4\lambda/(f\tau)+2$.

A spinor field $\psi$ on a Sasakian manifold $(\tilde M,\tilde \xi, \tilde g)$ 
with Reeb vector field $\tilde \xi$ is called a Sasakian quasi-Killing spinor of 
type $(a, b)$ if it  satisfies $\nabla_X\psi=aX\cdot\psi$ for 
$X\in\tilde\xi^\perp$ and $\nabla_{\tilde\xi}\psi=(a+b)\tilde\xi\cdot\psi$ for 
$a,b\in\R$. If $(\tilde M, \tilde \xi, \tilde g)$ admits a Sasakian 
quasi-Killing spinors of type $(a, b)$, then the scalar curvature $\S$ is 
constant and given by $\S=8m(2m+1)a^2+16mab$, see {\rm \cite{FriedKim00}}, Lemma 
6.4. In the following sense, in three dimensions, the converse is true. Let 
$(\tilde M,\tilde g, \tilde \xi)$  be a simply-connected three-dimensional 
Sasakian spin manifold with constant scalar curvature $\S$. Then there exist two 
linear independent Sasakian quasi-Killing spinors of type $(-1/2,3/4-\S/8)$, see 
{\rm \cite{FriedKim00}}, Theorem 8.4.

 We identify the spinor bundle $\tilde\Sigma N$ of $(N,\tilde g)$ with $\Sigma 
N$ by 
$\Sigma N \rightarrow\tilde\Sigma N$, $\ph\mapsto \tilde \ph$ such that a 
section $\ph$ in  $\Sigma N$ satisfies 
$$(X\clN \ph)^{\sim}=\tilde X\clN\tilde\ph,\quad 
(\nabla_X^N\ph)^\sim=\tilde\nabla^N_X\tilde\ph,$$
where $\tilde X:= X/|\tau|$ and $\tilde\nabla^N$ denotes the Levi-Civita 
connection on $\tilde \Sigma N$. Now we consider the restriction of the skew 
Killing spinor $\psi=\psi^++\psi^-$ to $N$. We will write $\psi^\pm$ instead of 
$\phi(\psi^\pm)$.  Then, by (\ref{Eres}), 
\begin{eqnarray*}
\tilde \nabla^N_\eta \tilde\psi^\pm&=&(\nabla^N_\eta \tilde\psi^\pm)^\sim 
=(-\frac\lambda{2f}\eta\clN\psi^{\pm})^\sim=-\frac\lambda{2f|\tau|}
\eta\clN(\psi^{\pm})^\sim, \\
\tilde \nabla^N_Z\tilde\psi^\pm&=&(\nabla^N_Z\tilde\psi^\pm)^\sim=(-\frac{\mu 
f}2Z\clN\psi^\pm)^\sim=-\frac{\mu 
f}{2|\tau|}Z\clN(\psi^\pm)^\sim=-\frac{\sgn(\tau)} 2 Z\clN(\psi^\pm)^\sim
\end{eqnarray*}
for $Z\in\eta^\perp$. Hence, $\tilde\psi^\pm$ is a Sasakian quasi-Killing spinor 
with $a=-\sgn(\tau)/2$ and 
$$b=-\frac\lambda{2f|\tau|}+\frac{\sgn(\tau)}2=\sgn(\tau)\Big(-\frac\lambda{
2f\tau}+\frac12\Big)=\sgn(\tau)\Big(\frac34-\frac \S8\Big).$$
Thus we are up to a change of orientation exactly in the situation described 
above. 

In dimension three Sasakian quasi-Killing spinors of type can also be understood 
as transversal Killing spinors, see {\rm\cite{GH08}} for a definition. If we 
return to our original metric $g$ on $N$, this means that the restrictions of 
$\psi^\pm$ to $N$ are transversal Killing spinors. Indeed, 
\begin{eqnarray*}
\bar\nabla_{\eta}\psi^\pm&=& \nabla^N_\eta\psi^\pm-\frac12\tau|\eta| s_3s_4\clN 
\psi^\pm\ =\ \nabla^N_\eta\psi^\pm+\frac12\tau|\eta| s_1\clN \psi^\pm\ =\ 
\Big(-\frac \lambda{2f}-\frac12\tau\Big)\eta\clN \psi^\pm,\\
\bar\nabla_{Z}\psi^\pm&=& \nabla^N_Z\psi^\pm-\frac12\tau s_1 J(Z)\clN \psi^\pm\ 
=\ -\frac{\mu f}2Z\clN\psi^\pm+\frac\tau 2 Z\clN \psi^\pm\ =\ 0
\end{eqnarray*}
holds for the transversal covariant derivative $\bar \nabla$ on $N$.
\end{erem}
\section*{Appendix: Doubly warped products}\label{s:dwp}

\begin{edefi}\label{d:dwp}
A \emph{doubly warped product} is a Riemannian manifold $(M,g)$ of the form 
$$(I\times \hat 
M,dt^2\oplus\rho(t)^2\hat{g}_{\hat{\eta}}\oplus\sigma(t)^2\hat{g}_{\hat{\eta}
^\perp}),$$ where  $(\hat M,\hat{g})$ is a Riemannian manifold with unit Killing 
vector field $\hat{\eta}$ and $\hat{g}_{\hat{\eta}}$, 
$\hat{g}_{\hat{\eta}^\perp}$ are the components of the metric $\hat{g}$ along 
$\R\hat{\eta}$ and $\hat{\eta}^\perp$, respectively, $I\subset\R$ is an open 
interval and $\rho,\sigma\colon I\to\R$ are smooth positive functions on 
$I$.\end{edefi}
\begin{edefi}\label{p:chardwpdistr}
Let $(M,g)$ be a Riemannian manifold.
A local DWP-structure $(\nu,\hat\eta)$ on $(M,g)$ consists of 
\begin{enumerate}
\item a unit geodesic vector field $\nu$ whose orthogonal complement 
distribution is integrable,
\item a nontrivial Killing vector field $\hat{\eta}$ on $(M,g)$ that is 
pointwise orthogonal to $\nu$ and whose length is constant along any integral 
leaf of $\nu^\perp$
\end{enumerate}
with the property that the Weingarten map $W:=-\nabla\nu$ of each integral leaf 
of $\nu^\perp$ has two eigenspaces, $\R\hat{\eta}$ and 
$\hat{\eta}^\perp\cap\nu^\perp$ and the corresponding eigenvalues $\lambda$ and 
$\mu$ are constant along the leaf. 
\end{edefi}

\begin{prop}
If $(M,g)$ is isometric to a doubly warped product, then $(M,g)$ admits a local 
DWP-structure. 
Conversely, if $(M,g)$ has a local DWP-structure, then it is locally isometric 
to a doubly-warped product.
\end{prop}

{\it Proof:} First assume that $(M,g)$ is isometric to a doubly warped product, 
thus $(M,g)=(I\times \hat 
M,dt^2\oplus\rho(t)^2\hat{g}_{\hat{\eta}}\oplus\sigma(t)^2\hat{g}_{\hat{\eta}
^\perp})$. Then we have the following expressions for the Levi-Civita connection 
$\nabla$ of $(M,g)$, see e.g. \cite[Sec. 3]{GinSemm11} (mind that our 
$\hat{\eta}$ here corresponds to $\hat{\xi}$ in \cite{GinSemm11} and that our 
$\rho$ and $\sigma$ correspond to $\rho\sigma$ and $\rho$, respectively).  For 
all sections $X,Y$ of $\pi_2^*Q$, where $Q:=\hat\eta^\perp\to \hat M$,
\begin{equation}\label{dwp}
{ \everymath={\displaystyle}
\begin{array}{lll} \displaystyle
\nabla_{\partial_t}\partial_t=0, 
&\nabla_{\partial_t}\hat{\eta}=\frac{\rho'}{\rho}\hat{\eta}, &\nabla_{\partial 
t}X=\frac{\partial X}{\partial t}+\frac{\sigma'}{\sigma}X, \\
\nabla_{\hat{\eta}}\partial_t=\frac{\rho'}{\rho}\hat{\eta}, 
&\nabla_{\hat{\eta}}\hat{\eta}=-\rho\rho'\partial_t, &
\nabla_{\hat{\eta}} 
X=\hat{\nabla}_{\hat{\eta}}X+\frac{\rho^2}{\sigma^2}\hat{h}X,\\
\nabla_X\partial_t=\frac{\sigma'}{\sigma}X, \ 
&\nabla_X\hat{\eta}=\frac{\rho^2}{\sigma^2}\hat{h}X, \ 
&\nabla_XY=\hat{\nabla}_XY-\frac{1}{\sigma^2}g(\hat{h}X,Y)\xi-\sigma'g(X,
Y)\partial_t,
\end{array}
}
\end{equation}
where $\hat{h}:=\nabla^{\hat{M}}\hat{\eta}\in\Gamma(\mathrm{End}(Q))$.
It is straightforward to see that $\nu:=\partial_t$ is a geodesic vector field 
with $\nu^\perp=T\hat M$, the vector field $\hat{\eta}$ (seen as a section of 
$\pi_2^*T\hat M\subset TM$) is Killing on $(M,g)$ with constant length along 
each $\{t\}\times \hat M$ and that 
$W:=-\nabla\partial_t=-\frac{\rho'}{\rho}\mathrm{Id}_{\R\hat{\eta}}\oplus-\frac{
\sigma'}{\sigma}\mathrm{Id}_Q$.

Conversely, let $(\nu,\hat{\eta})$ be a local DWP-structure on $(M,g)$. Let $p$ 
be a point in $M$. Then we find a local leaf $\hat M$ of $\nu^\perp$ such that 
the integral curves of $\nu$ starting from $\hat M$ are defined at least on an 
interval $(-t_0,t_0)$. We denote by $\hat{g}$ the induced metric on $\hat M$. Up 
to rescaling $\hat{\eta}$ by a nonzero constant, we may assume that 
$\hat{g}(\hat{\eta},\hat{\eta})=1$ along $\hat M$.
Consider the map $F\colon (-t_0,t_0)\times \hat M\to M$ given by 
$F(t,x):=F_t(x)$, where $(F_t)_t$ is the flow of the vector field $\nu$. The map 
$F$ is clearly a local diffeomorphism. Next we identify the pull-back metric 
$F^*g$ on $(-t_0,t_0)\times \hat M$. For any given $(t,x)\in(-t_0,t_0)\times 
\hat M$ and $X\in T_x \hat M$, we have 
\[(F^* 
g)_{(t,x)}(\partial_t,X)=g_{F(t,x)}(\nu,d_xF_t(X))=g_{F(t,x)}(d_xF_t(\nu),
d_xF_t(X))=(F_t^*g)_x(\nu,X).\]
Since $\nu$ is geodesic of constant length, 
$(\mathcal{L}_{\nu}g)(\nu,Y)=g(\nabla_\nu\nu,Y)+g(\nabla_Y\nu,\nu)=0$ holds for 
all  $Y\in TM$. Consequently, the derivative
$$ 
\frac{\partial}{\partial t}(F_t^*g)_x(\nu_,X)=\frac{\partial}{\partial 
s}(F_{t+s}^*g)_x(\nu_,X)\big|_{s=0}=(\mathcal{L}_{\nu}g)_{F(t,x)}((F_t)_*\nu,
(F_t)_*X)=(\mathcal{L}_{\nu}g)_{F(t,x)}(\nu,(F_t)_*X)
$$
vanishes, thus $(F_t^*g)_x(\nu,X)=(F_0^*g)_x(\nu,X)=g_x(\nu,X)=0$ for all 
$(t,x)\in (-t_0,t_0)\times \hat M$.
This proves the splitting $F^*g=dt^2\oplus g_t$, where $g_t:=(F_t^*g){|_{T\hat 
M\times T\hat M}}$.
As a next step, we compute $g_t$ more precisely along each of the distributions 
$\R\hat{\eta}$ and $Q$ of $T\hat M$.
We first notice that $\hat{\eta}$ is invariant under the flow of $\nu$.
Namely, we write $W=\lambda\mathrm{Id}_{\R\eta}\oplus\mu\mathrm{Id}_Q$ for 
functions $\lambda,\mu\colon\, \R\to\R$, which are constant along each integral 
leaf of $\nu^\perp$ by assumption. Since $\hat{\eta}$ is Killing, we have 
$g(\nabla_\nu\hat{\eta},\nu)=0$.
Moreover, because of $\hat{\eta}\perp\nu$,
\[g(\nabla_\nu\hat{\eta},\hat{\eta})=-g(\nabla_{\hat{\eta}}\hat{\eta},
\nu)=g(\nabla_{\hat{\eta}}\nu,\hat{\eta})=-g(W\hat{\eta},\hat{\eta})=-\lambda 
g(\hat{\eta},\hat{\eta}).\]
Note that this proves in particular that, if $\hat{\eta}$ vanishes at a point, 
then it must vanish on the corresponding integral leaf of $\nu^\perp$ and 
therefore identically on the image of $F$ since $g(\hat{\eta},\hat{\eta})$ 
satisfies the ODE $\nu(g(\hat{\eta},\hat{\eta}))=-2\lambda 
g(\hat{\eta},\hat{\eta})$. Furthermore, for every $X\in Q$,
\[g(\nabla_\nu\hat{\eta},X)=-g(\nabla_X\hat{\eta},\nu)=g(\nabla_X\nu,\hat{\eta}
)=-\mu g(\hat{\eta},X)=0.\]
As a first consequence, $\nabla_\nu\hat{\eta}=-\lambda\hat{\eta}$. This implies
$\mathcal{L}_{\nu}\hat \eta= 
[\nu,\hat{\eta}]=\nabla_\nu\hat{\eta}-\nabla_{\hat{\eta}}\nu=0,$
so that 
$(F_t)_*\hat{\eta}=\hat{\eta}$ for every $t\in \R$. For any $X,Y\in \nu^\perp$, 
we have 
\be 
(\mathcal{L}_\nu g)(X,Y)&=&g(\nabla_X\nu,Y)+g(\nabla_Y\nu,X)=-2g(WX,Y).
\ee
Thus, in particular,
\begin{eqnarray} 
\frac{\partial}{\partial 
s}(F_s^*g)(\hat\eta,Y){\big|_{s=t}}&=&(\mathcal{L}_{\nu}g)_{F(t,x)}
((F_t)_*\hat\eta,(F_t)_*Y)=(\mathcal{L}_{\nu}g)_{F(t,x)}(\hat\eta,
(F_t)_*Y)\nonumber\\
&=&-2g(W\hat\eta,(F_t)_*Y) =-2(\lambda\circ F)\cdot g(\hat \eta, (F_t)_*Y) 
\label{final}\\
&=&-2(\lambda\circ F) (F_t^*g)(\hat \eta, Y). \nonumber
\end{eqnarray}
Consequently, for fixed $Y\in T_x\hat M\cap \nu^\perp$, the function 
$\varphi(t):=(F_t)^*g(\hat \eta, Y)$ satisfies the differential equation 
$$\varphi'(t)=-2\lambda(F(t,x))\cdot\varphi(t).$$
For $Y\in Q$, we have $\varphi(0)=0$, thus $\varphi=0$. This means that the flow 
of $\nu$ preserves the distribution~$\hat\eta^\perp$. For $Y=\hat \eta$, we have 
$\varphi(0)= 1$, thus 
$$(F_t)^*g(\hat \eta, \hat \eta)=\varphi(t)=\exp(-2\int_0^t\lambda\circ F_s 
ds)=:\rho(t)^2.$$
Finally, for $X,Y\in Q$, a computation analogous to (\ref{final}) shows that 
\[\frac{\partial}{\partial s}(F_s^*g)(X,Y){\Big|_{s=t}}=-2(\mu\circ 
F)(F_t^*g)(X,Y),\]
which yields $$(F_t^*g)(X,Y)=\exp(-2\int_0^t\mu\circ F_s ds)\cdot \hat 
g(X,Y)=:\sigma(t)^2\hat g(X,Y).$$
It remains to notice that $\hat{\eta}$ must be a Killing vector field along 
$(\hat  M,\hat{g})$ since it is already Killing on $(M,g)$ and is tangent to 
$\hat M$.
On the whole, we obtain the doubly warped product metric as required.
\findemo

\end{document}